\documentclass[12pt]{article}

\usepackage{latexsym}
\usepackage{amssymb}

\newtheorem{thm}{Theorem}
\newtheorem{lem}{Lemma}
\newtheorem{prop}{Proposition}

\newtheorem{dfn}{Definition}

\newtheorem{example}{Example}

\begin{document}
{

\begin{center}
{\large\bf
ON THE TRUNCATED MATRIX POWER MOMENT PROBLEM WITH AN OPEN GAP}
\end{center}

\begin{center}
S. M. Zagorodnyuk
\end{center}

\section{Introduction.}

The main objects of our present investigation are some truncated matrix power moment problems. At first, we study \textit{the truncated matrix Hamburger
moment problem with an odd number of moments}: find a left-continuous non-decreasing matrix function $M(x) = ( m_{k,l}(x) )_{k,l=0}^{N-1}$
on $\mathbb{R}$, $M(-\infty)=0$, such that
\begin{equation}
\label{f1_5}
\int_\mathbb{R} x^n dM(x) = S_n,\qquad n=0,1,\ldots,2d,
\end{equation}
where $\{ S_n \}_{n=0}^{2d}$ is a prescribed sequence of Hermitian $(N\times N)$ complex matrices (moments), $d\in \mathbb{Z}_+$, $N\in \mathbb{N}$.

\noindent
Secondly, we study the moment problem~(\ref{f1_5}) with an additional constraint posed on the matrix measure $M(\delta)$ (generated by the
distribution function $M(x)$):
\begin{equation}
\label{f1_10}
M(\Delta) = 0,
\end{equation}
where $\Delta$ is a given open subset of $\mathbb{R}$ (a gap).

\noindent
Observe that the moment problem with a gap~(\ref{f1_5}),(\ref{f1_10}) in various particular cases coincides with the well-known moment problems.
For instance, the case $\Delta=(-\infty,0)$ leads to the truncated matrix Stieltjes moment problem; the case
$\Delta=(-\infty,a)\cup (b,+\infty)$, $-\infty < a < b < +\infty$, leads to the truncated matrix Hausdorff moment problem, etc.

Moment problems form a particular and very important case of interpolation problems for various classes of analytic functions.
For the theory of classical scalar moment problems see, e.g., books~\cite{cit_600_Akh}, \cite{cit_1000_Ber}.
Matrix moment problems appeared as a natural generalizations of classical problems.
The history for the truncated matrix Hamburger moment problem~(\ref{f1_5}) was presented in our paper~\cite{cit_35000_Z}.
That paper will be intensively used in the present paper. The results of paper~\cite{cit_35000_Z} on the truncated matrix Hamburger moment problem
will be our starting point here. Regarding the history for the truncated matrix Stieltjes and Hausdorff moment problems we refer to
papers~\cite{cit_2000_D}, \cite{cit_300_A_T}, \cite{cit_3000_H_C}, \cite{cit_37000_Z} and papers cited there.
For scalar truncated power moment problems with intervals as gaps see~\cite{cit_1500_D_M} and references therein.

The moment problem~(\ref{f1_5}) (the moment problem~(\ref{f1_5}),(\ref{f1_10})) is said to be \textbf{determinate} if it has a unique solution,
and \textbf{indeterminate} if it has more than one solution.
In Section~2 we obtain some necessary and sufficient conditions for the determinacy of the truncated matrix Hamburger moment problem~(\ref{f1_5}).
Then we derive a Nevanlinna-type parametrization for all solutions of the moment problem~(\ref{f1_5}).
The main tools are properties of the corresponding generalized resolvents.
In the case of a determinate moment problem~(\ref{f1_5}) the unique solution is constructed explicitly, as well.

In Section~3 we give a criterion for the solvability of the truncated matrix power moment problem with a gap~(\ref{f1_5}),(\ref{f1_10}).
A Nevanlinna-type formula for all solutions of the moment problem is obtained, as well.

We emphasize that the coefficients of linear fractional transformations in all of the presented Nevanlinna-type formulas can be
calculated by the prescribed moments explicitly.

\noindent
\textbf{Notations.}  As usual, we denote by $\mathbb{R}, \mathbb{C}, \mathbb{N}, \mathbb{Z}, \mathbb{Z}_+$
the sets of real numbers, complex numbers, natural numbers, integers and non-negative integers,
respectively; $\mathbb{R}_e = \mathbb{C}\backslash \mathbb{R}$; $\mathbb{C}_\pm = \{ z\in \mathbb{C}:\ \pm\mathop{\rm Im}\nolimits z > 0\}$;
$\mathbb{T} = \{ z\in \mathbb{C}:\ |z|=1 \}$.
By $\Pi_z$ we mean the set, either $\mathbb{C}_+$ or $\mathbb{C}_-$, which contains $z$; $z\in \mathbb{R}_e$.
We set $\Pi_{\lambda}^\varepsilon = \{ z\in\Pi_{\lambda}:\ \varepsilon < |\arg z| < \pi - \varepsilon \},\qquad
0 < \varepsilon < \frac{\pi}{2}$, $\lambda\in \mathbb{R}_e$.
By $\mathfrak{B}(\mathbb{R})$ we mean a set of all Borel subsets of $\mathbb{R}$.

Let $N,K\in \mathbb{N}$.
The set of all complex matrices of order $(N\times K)$ we denote by~$\mathbb{C}_{N\times K}$.
The set of Hermitian non-negative complex matrices of order $(N\times N)$ will be denoted by~$\mathbb{C}_{N\times N}^{\geq}$.
If $M\in \mathbb{C}_{N\times K}$ then $M^*$ means the complex conjugate matrix.
If $L\in \mathbb{C}_{N\times N}$ then $\mathop{\rm Ker}\nolimits L = \{ x\in \mathbb{C}_{N\times 1}:\ Lx=0 \}$.
By $I_N$ we denote the identity matrix of size $(N\times N)$.
By $\mathbb{C}^N$ ($\mathbb{C}_N$) we denote the complex Euclidean space of vectors from $\mathbb{C}_{1\times N}$ (respectively
from $\mathbb{C}_{N\times 1}$).
We denote $\vec e_n = (\delta_{n,0},\delta_{n,1},\ldots,\delta_{n,N-1})\in \mathbb{C}_{1\times N}$, $0\leq n\leq N-1$.

\noindent
Let $M(x)$ be a left-continuous non-decreasing matrix function $M(x) = ( m_{k,l}(x) )_{k,l=0}^{N-1}$
on $\mathbb{R}$, $M(-\infty)=0$, and $\tau_M (x) := \sum_{k=0}^{N-1} m_{k,k} (x)$;
$M'_\tau (x) = ( dm_{k,l}/ d\tau_M )_{k,l=0}^{N-1}$.  We denote by $L^2(M)$ a set (of classes of equivalence)
of vector functions $f: \mathbb{R}\rightarrow \mathbb{C}_{1\times N}$, $f = (f_0,f_1,\ldots,f_{N-1})$, such that
$$ \| f \|^2_{ L^2(M) } := \int_\mathbb{R}  f(x) M'_\tau (x) f^*(x) d\tau_M (x) < \infty. $$
The space $L^2(M)$ is a Hilbert space with the scalar product
$$ ( f,g )_{L^2(M)} := \int_\mathbb{R}  f(x) M'_\tau (x) g^*(x) d\tau_M (x),\qquad f,g\in L^2(M). $$

For a separable Hilbert space $H$ we denote by $(\cdot,\cdot)_H$ and $\| \cdot \|_H$ the scalar
product and the norm in $H$, respectively. In obvious cases, the indices may be omitted.

\noindent
For a linear operator $\mathcal{A}$ in $H$ we denote by $D(\mathcal{A})$ its domain, by $R(\mathcal{A})$ its range, and by
$\mathcal{A}^*$ we denote its adjoint if it exists. If $\mathcal{A}$ is bounded, then $\| \mathcal{A} \|$ stands for its operator norm.
For a set $S$ in $H$, we
denote by $\mathop{\rm Lin}\nolimits S$ and $\mathop{\rm span}\nolimits S$ the linear span and the closed
linear span, in the norm of $H$, of $S$, respectively.
For a set $M\subseteq H$ we denote by $\overline{M}$ the closure of $M$ with respect to the norm of $H$.
By $E_H$ we denote the identity operator in $H$, i.e. $E_H x = x$, $x\in H$.
If $H_1$ is a subspace of $H$, by $P_{H_1} = P_{H_1}^{H}$ we denote the operator of the orthogonal projection on $H_1$
in $H$.

By $\mathcal{S}(D;N,N')$ we denote a class of all analytic in a domain $D\subseteq \mathbb{C}$
operator-valued functions $F(z)$, which values are linear non-expanding operators mapping the whole
$N$ into $N'$, where $N$ and $N'$ are some Hilbert spaces.

Consider a closed symmetric operator $A$ in a Hilbert space $H$.
Denote
$\mathcal{M}_z = \mathcal{M}_z(A) = (A - zE_H) D(A)$,
$\mathcal{N}_z = \mathcal{N}_z(A) = H\ominus \mathcal{M}_z(A)$, where $z\in \mathbb{C}$.

Choose and fix an arbitrary point $\lambda_0\in \mathbb{R}_e$.
Set $U_{\lambda_0}(A) = (A - \overline{\lambda_0}E_H)(A - \lambda_0 E_H)^{-1}$.
Define a linear operator $X_{\lambda_0}=X_{\lambda_0}(A)$ in the following way:
\begin{equation}
\label{f1_5_p2_1}
D(X_{\lambda_0}) = P^H_{\mathcal{N}_{\lambda_0}(A)} (H\ominus \overline{D(A)}),
\end{equation}
\begin{equation}
\label{f1_6_p2_1}
X_{\lambda_0} P^H_{\mathcal{N}_{\lambda_0}(A)} h = P^H_{\mathcal{N}_{\overline {\lambda_0}}(A)} h,\qquad h\in H\ominus \overline{D(A)}.
\end{equation}
The operator $X_{\lambda_0}=X_{\lambda_0}(A)$ is called \textit{forbidden with respect to $A$}.

A function $F(\lambda)\in \mathcal{S}(\Pi_{\lambda_0}; \mathcal{N}_{\lambda_0}(A),
\mathcal{N}_{ \overline{\lambda_0} }(A))$ is said to be  {\it $\lambda_0$-admissible (admissible) with respect to the operator $A$},
if the validity of
\begin{equation}
\label{f11_1_p2_1}
\lim_{\lambda\in\Pi_{\lambda_0}^\varepsilon,\ \lambda\to\infty} F(\lambda) \psi = X_{\lambda_0}\psi,
\end{equation}
\begin{equation}
\label{f11_2_p2_1}
\underline{\lim}_{\lambda\in\Pi_{\lambda_0}^\varepsilon,\ \lambda\to\infty}
\left[
|\lambda| (\| \psi \|_H - \| F(\lambda) \psi \|_H)
\right] < +\infty,
\end{equation}
for some $\varepsilon$: $0<\varepsilon <\frac{\pi}{2}$,
implies $\psi = 0$.

A set of all operator-valued functions $F(\lambda)\in \mathcal{S}(\Pi_{\lambda_0}; \mathcal{N}_{\lambda_0}(A),
\mathcal{N}_{ \overline{\lambda_0} }(A))$, which are $\lambda_0$-admissible with respect to the operator $A$,
we shall denote by
$$ \mathcal{S}_{a; \lambda_0} (\Pi_{\lambda_0}; \mathcal{N}_{\lambda_0}(A),
\mathcal{N}_{ \overline{\lambda_0} }(A)) = \mathcal{S}_{a} (\Pi_{\lambda_0}; \mathcal{N}_{\lambda_0}(A),
\mathcal{N}_{ \overline{\lambda_0} }(A)). $$

For a closed isometric operator $V$ in a Hilbert space $H$ we denote:
$M_\zeta(V) = (E_H - \zeta V) D(V)$, $N_\zeta(V) = H\ominus M_\zeta(V)$, $\zeta\in \mathbb{C}$; $M_\infty(V)=R(V)$, $N_\infty(V)= H\ominus R(V)$.

\section{The truncated matrix Hamburger moment problem with an odd number of moments.}
Consider the moment problem~(\ref{f1_5}) with a prescribed sequence $\{ S_n \}_{n=0}^{2d}$ of Hermitian $(N\times N)$ complex matrices,
$d,N\in \mathbb{N}$.
Set
\begin{equation}
\label{f2_3}
\Gamma_d = \left( S_{i+j} \right)_{i,j=0}^d,\qquad \widehat\Gamma_{d-1} = (S_{i+j+2})_{i,j=0}^{d-1}.
\end{equation}
The following conditions:
\begin{equation}
\label{f2_7}
\Gamma_{d}\geq 0,\quad \mathop{\rm Ker}\nolimits \Gamma_{d-1} \subseteq \mathop{\rm Ker}\nolimits \widehat\Gamma_{d-1},
\end{equation}
are necessary and sufficient for the solvability of the moment problem~(\ref{f1_5}), see~Remark on page~286 in~\cite{cit_35000_Z}.

\noindent
Suppose that conditions~(\ref{f2_7}) are satisfied.
Let $\Gamma_d = (\gamma_{n,m}^d)_{n,m=0}^{dN +N-1}$, $\gamma_{n,m}=\gamma_{n,m}^d\in \mathbb{C}$, be the usual representation of the block
matrix $\Gamma_d$.
Repeating a construction from~\cite[p. 281]{cit_35000_Z}
we get a finite-dimensional Hilbert space $H$ and
a sequence $\{ x_n \}_{n=0}^{dN +N-1}$ in $H$, such that
\begin{equation}
\label{f2_10}
(x_n,x_m)_H = \gamma_{n,m},\qquad n,m=0,1,...,dN +N-1,
\end{equation}
and $\mathop{\rm span}\nolimits \{ x_n \}_{n=0}^{dN +N-1} = H$.

\noindent
Set $L_a = \mathop{\rm Lin}\nolimits \{ x_n \}_{n=0}^{dN-1}$,
and consider the following linear operator $A$ with the domain $D(A)=L_a$:
\begin{equation}
\label{f2_15}
A x = \sum_{k=0}^{dN-1} \alpha_k x_{k+N},\qquad x\in L_a,\ x = \sum_{k=0}^{dN-1} \alpha_k x_{k},\ \alpha_k\in \mathbb{C}.
\end{equation}
In particular, we have
\begin{equation}
\label{f2_20}
A x_k =  x_{k+N},\qquad 0\leq k\leq dN-1.
\end{equation}
In~\cite{cit_35000_Z} it is checked that the operator $A$ is well-defined and symmetric.
Since it acts in a finite-dimensional space, it is closed and its defect numbers are equal.
By Theorem~5 in~\cite{cit_35000_Z} all solutions of the moment problem~(\ref{f1_5}) have the following form:
\begin{equation}
\label{f2_25}
M(\lambda) = (m_{k,j} (\lambda))_{k,j=0}^{N-1},\quad
m_{k,j} (\lambda) = ( \mathbf E_\lambda x_k, x_j)_H,
\end{equation}
where $\mathbf E_\lambda$ is a left-continuous spectral function of the operator $A$.
Moreover, the correspondence between all left-continuous spectral functions of $A$ and all solutions
of the moment problem is one-to-one.

\noindent
Since there exists a one-to-one correspondence between all
generalized resolvents of $A$  and all left-continuous spectral functions:
$$
(\mathbf R_z f,g)_H = \int_\mathbb{R} \frac{1}{\lambda - z} d( \mathbf E_\lambda f,g)_H,\qquad f,g\in H,\ z\in \mathbb{R}_e,
$$
we conclude that
all solutions of the moment problem~(\ref{f1_5}) have the following form:
\begin{equation}
\label{f2_30}
M(\lambda) = (m_{k,j} (\lambda))_{k,j=0}^{N-1},\quad
\int_\mathbb{R} \frac{1}{\lambda - z} dm_{k,j} (\lambda) = ( \mathbf R_z x_k, x_j)_H,\quad z\in \mathbb{R}_e,
\end{equation}
where $\mathbf R_z$ is a generalized resolvent of $A$. On the other hand, for an arbitrary generalized resolvent there corresponds by~(\ref{f2_30})
a solution $M(\lambda)$ of the moment problem (via the Stieltjes-Perron inversion formula). Moreover, for different
generalized resolvents there correspond by~(\ref{f2_30})
different solutions of the moment problem.%ok

Let us apply the Gram-Schmidt orthogonalization procedure to the following sequence:
$$ x_0, x_1, ..., x_{dN-1}, x_{dN}, ..., x_{dN+N-1}, $$
removing the linear dependent elements if they appear.
During the orthogonalization of the first $dN$ elements ($\{ x_k \}_{k=0}^{dN-1}$) we shall obtain $\kappa$ orthonormal elements
$\mathfrak{A}_\infty := \{ f_j \}_{j=0}^{\kappa-1}$, $0\leq\kappa\leq dN$. The case $\kappa = 0$ means that $x_k=0$, $0\leq k\leq dN-1$,
and $\mathfrak{A}_\infty = \emptyset$. By~(\ref{f2_10}) we see  that in this case holds $S_0 = 0$, and therefore $M(x)\equiv 0$ and
the moment problem is determinate. In this case all given moments $S_n$ are zero matrices. Then by~(\ref{f2_10}) it follows that
all $x_k$ are zero elements and $H = \{ 0 \}$.

Applying the orthogonalization to the rest of the elements ($\{ x_k \}_{k=dN}^{dN+N-1}$) we shall obtain $\kappa'$
orthonormal elements
$\mathfrak{A}_\infty' := \{ f_j' \}_{j=0}^{\kappa'-1}$, $0\leq\kappa'\leq N$.
In the case $\kappa'=0$ we mean $\mathfrak{A}_\infty' = \emptyset$.

Notice that $\mathfrak{A}_\infty$ ($\mathfrak{A}_\infty'$) is an orthonormal basis in $D(A)$ (respectively, in $H\ominus D(A)$), if
it is non-empty.

Observe that the elements of $\mathfrak{A}_\infty$ and $\mathfrak{A}_\infty'$ are \textbf{linear combinations of $x_j$ with some coefficients,
which can be explicitly calculated by~(\ref{f2_10})}.
Similar features will be true for all orthogonalization procedures in this paper.

\begin{thm}
\label{t2_1}
Let the truncated matrix Hamburger moment problem~(\ref{f1_5}) be given and conditions~(\ref{f2_7}) be satisfied.
Let the operator $A$ in the Hilbert space $H$ be constructed as in~(\ref{f2_15}).
The following conditions are equivalent:

\begin{itemize}

\item[{\rm (A)}]
The moment problem~(\ref{f1_5}) is determinate;

\item[{\rm (B)}]
The operator $A$ is self-adjoint;

\item[{\rm (C)}] $\mathfrak{A}_\infty' = \emptyset$.

\end{itemize}

If the above conditions are satisfied then the unique solution of the moment problem~(\ref{f1_5})
is given by the relation~(\ref{f2_25}) where $\mathbf{E}_\lambda$ is the orthogonal left-continuous spectral function of the self-adjoint operator $A$.
This solution is piecewise constant.

\end{thm}
\textbf{Proof.}
(A)$\Rightarrow$(B).
Suppose to the contrary that the moment problem~(\ref{f1_5}) is determinate and the operator $A$ has non-zero defect numbers.
We can choose non-zero elements $u_1\in N_i(A)$ and $u_2\in N_{-i}(A)$, $\| u_1 \|_H = \| u_2 \|_H = 1$,
and set
$$ F_1(\lambda) = 0,\quad F_2(\lambda) = \frac{1}{2} (\cdot,u_1)_H u_2,\qquad \lambda\in \mathbb{C}_+. $$
The operator-valued functions $F_1(\lambda)$ and $F_2(\lambda)$ %are admissible with respect to $A$ and
generate by Theorem~7 in~\cite{cit_35000_Z}
different solutions of the moment problem~(\ref{f1_5}).
We obtained a contradiction.

\noindent
(B)$\Rightarrow$(A). The self-adjoint operator $A$ has a unique (left-continuous) spectral function, therefore by~(\ref{f2_25})
we conclude that the moment problem is determinate.

\noindent
(B)$\Rightarrow$(C). Since $A$ is self-adjoint and it acts in a finite-dimensional space, then $D(A)=H$.
Then $\mathfrak{A}_\infty' = \emptyset$.

\noindent
(C)$\Rightarrow$(B). If $\mathfrak{A}_\infty' = \emptyset$, then $D(A)=H$. Therefore $A$ is self-adjoint.

\noindent
$\Box$

\begin{example}
\label{e2_1}
Consider the moment problem~(\ref{f1_5}) with $d=1$, $N=2$, and the following moments:
$$ S_0 =
\left(
\begin{array}{cc} \frac{1}{3} & 0 \\
0 & 1\end{array}
\right),\quad
S_1 =
\left(
\begin{array}{cc} \frac{1}{2} & 0 \\
0 & 1\end{array}
\right),\quad
S_2 =
\left(
\begin{array}{cc} 1 & 0 \\
0 & 1\end{array}
\right). $$
Conditions~(\ref{f2_7}) are verified directly.
Consider a Hilbert space $H$ and a sequence $\{ x_n \}_{n=0}^3$ in $H$ such that relation~(\ref{f2_10}) holds and
$\mathop{\rm Lin}\nolimits \{ x_n \}_{n=0}^{3} = H$.
Let us apply the Gram-Schmidt orthogonalization procedure to the following sequence:
$$ x_0, x_1, x_2, x_3, $$
removing the linear dependent elements if they appear.

\noindent
\textbf{Step~1)} Calculate
$$ \| x_0 \|_H^2 = (x_0,x_0)_H = \gamma_{0,0} = \frac{1}{3}. $$
Then $f_0 := \sqrt{3} x_0$.

\noindent
\textbf{Step~2)} Calculate
$$ h:= x_1 - (x_1,f_0)_H f_0 = x_1 - 3\gamma_{1,0} x_0 = x_1; $$
$$ \| h \|_H^2 = (x_1,x_1)_H = \gamma_{1,1} = 1. $$
Then
$$ f_1 = \frac{h}{ \| h \|_H } = x_1. $$
We get $\mathfrak{A}_\infty = \{ f_0, f_1 \}$, $\kappa = 2$.

\noindent
\textbf{Step~3)} Calculate
$$ \widehat h:= x_2 - (x_2,f_0)_H f_0 - (x_2,f_1)_H f_1 = x_2 - \frac{3}{2} x_0; $$
$$ \| \widehat h \|_H^2 = \left( x_2 - \frac{3}{2} x_0, x_2 - \frac{3}{2} x_0 \right)_H = \frac{1}{4}. $$
Then
$$ f_0' = \frac{ \widehat h }{ \| \widehat h \|_H } = 2x_2 - 3x_0. $$

\noindent
\textbf{Step~4)} Calculate
$$ \widetilde h:= x_3 - (x_3,f_0)_H f_0 - (x_3,f_1)_H f_1 - (x_3,f_0')_H f_0' = x_3 - x_1; $$
$$ \| \widetilde h \|_H^2 = \left( x_3 - x_1, x_3 - x_1 \right)_H = 0. $$
Thus, $x_3\in \mathop{\rm Lin}\nolimits \{ f_0, f_1, f_0' \}$.
Consequently, we get $\mathfrak{A}_\infty' = \{ f_0' \}$, $\kappa' = 1$.

\noindent
By Theorem~\ref{t2_1} we obtain that the moment problem is indeterminate.

\end{example}

Consider the \textbf{indeterminate} moment problem~(\ref{f1_5}) with a prescribed sequence $\{ S_n \}_{n=0}^{2d}$ of Hermitian $(N\times N)$ complex matrices,
$d,N\in \mathbb{N}$.
Our aim is to derive a Nevanlinna-type parameterization for this moment problem. A general idea of this derivation is similar to the
idea used in our paper~\cite{cit_39000_Z} for the full matrix Hamburger moment problem. We shall need some auxiliary results for
the generalized resolvents, as well.

Set
$$ y_k := (A-iE_H) x_k = x_{k+N} - ix_k,\qquad k=0,1,...,dN-1; $$
$$ H^- := (A-iE_H) D(A) = \mathop{\rm Lin}\nolimits \{ y_k \}_{k=0}^{dN-1},\quad H^+ := (A+iE_H) D(A). $$
Let us apply the Gram-Schmidt orthogonalization procedure to the following sequence:
\begin{equation}
\label{f2_32}
y_0, y_1, ..., y_{dN-1}, x_{0}, x_1, ..., x_{N-1},
\end{equation}
removing the linear dependent elements if they appear. Observe that the linear span of the above elements coincide with $H$, as
it follows from the formula~(81) in~\cite[p. 284]{cit_35000_Z}.

\noindent
During the orthogonalization of the first $dN$ elements ($\{ y_k \}_{k=0}^{dN-1}$) we shall obtain $\tau$ orthonormal elements
$\mathfrak{A} := \{ u_j \}_{j=0}^{\tau-1}$, $1\leq\tau\leq dN$. In fact, if  we would have $y_k=0$, $0\leq k\leq dN-1$,
then $H^- = \{ 0 \}$ and $D(A)=\{ 0 \}$. Then $x_0 = x_1 = ... =x_{N-1} = 0$, and by~(\ref{f2_10}) we get $S_0 = 0$. Then
$M(x)\equiv 0$, and the moment problem is determinate, what contradicts to our assumptions.
Thus, at least one orthonormal element will be constructed during the orthogonalization.

Applying the orthogonalization to the rest of the elements ($\{ x_k \}_{k=0}^{N-1}$) we shall obtain $\delta$
orthonormal elements
$\mathfrak{A}' := \{ u_j' \}_{j=0}^{\delta-1}$, $1\leq\delta\leq N$. If we would not obtain new elements during this orthogonalization,
then $H^- = H$, and the operator $A$ would be self-adjoint. By Theorem~\ref{t2_1} we would obtain that the moment problem is determinate and
we would obtain a contradiction.

Observe that $\mathfrak{A}$ ($\mathfrak{A}'$) is an orthonormal basis in $H^-$ (respectively, in $H\ominus H^-$).
Consider the Cayley transformation of $A$:
$$ V := (A+iE_H)(A-iE_H)^{-1}. $$
Set
$$ v_k = V u_k,\qquad k=0,1,...,\tau - 1. $$
Notice that $\mathfrak{A}_v := \{ v_k \}_{k=0}^{\tau - 1}$ is an orthonormal basis in $H^+$.

Let us apply the Gram-Schmidt orthogonalization procedure to the following sequence:
$$ v_0, v_1, ..., v_{\tau - 1}, x_{0}, x_1, ..., x_{N-1}, $$
removing the linear dependent elements if they appear.
Notice that the linear span of the above elements coincide with $H$, as
it follows from the formula~(81) in~\cite[p. 284]{cit_35000_Z}.

Observe that the first $\tau$ elements are already orthonormal.
Applying the orthogonalization to the rest of the elements ($\{ x_k \}_{k=0}^{N-1}$) we shall obtain $\delta$
orthonormal elements
$\mathfrak{A}_v' := \{ v_k' \}_{k=0}^{\delta-1}$.
Notice that $\mathfrak{A}_v'$ is an orthonormal basis in $H\ominus H^+$.

We shall need the following description of the generalized resolvents of a symmetric operator.
Also, we think that this description has some interest in general, since it uses bounded operators and allows to calculate
the matrix of the generalized resolvent in an orthonormal basis.
\begin{prop}
\label{p2_1}
Let $A$ be a closed symmetric operator in a Hilbert space $H$, and $z\in \mathbb{R}_e$ be an arbitrary
point.
An arbitrary generalized resolvent $\mathbf{R}_{s;\lambda}$ of the operator $A$ has the following form:
$$
\mathbf R_{s;\lambda}(A) =
\frac{ z - \overline{z} }{ (\lambda - \overline{z}) (\lambda - z) }
\left(
\left[
E_H - \frac{\lambda - z}{\lambda - \overline{z}}
\left(
U_z(A) \oplus F(\lambda)
\right)
\right]^{-1}
-
\frac{\lambda - \overline{z}}{z - \overline{z}} E_H
\right), $$
\begin{equation}
\label{f2_40}
\lambda\in\Pi_z\backslash\{ z \},
\end{equation}
where $F(\lambda)$ is a function from $\mathcal{S}_{a;z}(\Pi_z; \mathcal{N}_z(A),
\mathcal{N}_{\overline{z}}(A))$.
Conversely, an arbitrary function $F(\lambda)\in \mathcal{S}_{a;z}(\Pi_{z}; \mathcal{N}_{z}(A),
\mathcal{N}_{\overline{z}}(A))$
defines by relation~(\ref{f2_40}) a generalized resolvent
$\mathbf{R}_{s;\lambda}$ of the operator $A$.
Moreover, for different functions from
$\mathcal{S}_{a;z}(\Pi_{z}; \mathcal{N}_{z}(A),
\mathcal{N}_{\overline{z}}(A))$
there correspond different generalized resolvents of the operator $A$.
\end{prop}
\textbf{Proof.}
By Theorem~2.9 in~\cite[p. 199]{cit_46000_Z} the following relation:
\begin{equation}
\label{f2_46}
\mathbf{R}_{u;\frac{\lambda - z}{\lambda - \overline{z}}} (U_z)
= \frac{\lambda - \overline{z}}{z-\overline{z}} E_H +
\frac{(\lambda - \overline{z})(\lambda - z)}{z-\overline{z}} \mathbf{R}_{s;\lambda} (A),\qquad
\lambda\in \mathbb{R}_e\backslash\{ z,\overline{z} \},
\end{equation}
establishes a one-to-one correspondence between all generalized resolvents $\mathbf{R}_{s;\lambda}(A)$
of the operator $A$
and those generalized resolvents $\mathbf{R}_{u;\zeta} (U_z)$ of the closed isometric
operator $U_z$ which are generated by extensions of $U_z$ without non-zero fixed points.
Let restrict relation~(\ref{f2_46}) to $\lambda\in \Pi_z\backslash\{ z \}$, and express
the generalized resolvent $\mathbf{R}_{s;\lambda} (A)$:
\begin{equation}
\label{f2_48}
\mathbf{R}_{s;\lambda} (A) =
\frac{ z - \overline{z} }{ (\lambda - \overline{z}) (\lambda - z) }
\left(
\mathbf{R}_{u;\frac{\lambda - z}{\lambda - \overline{z}}} (U_z)
-
\frac{\lambda - \overline{z}}{z - \overline{z}} E_H
\right),\quad \lambda\in\Pi_z\backslash\{ z \}.
\end{equation}
Let $F(\lambda)\in \mathcal{S}_{a;z}(\Pi_z; \mathcal{N}_z(A),
\mathcal{N}_{\overline{z}}(A))$ be the parameter in the Shtraus formula for the generalized resolvents
(see, e.g.~\cite[Theorem 3.34]{cit_46000_Z}) corresponding to $\mathbf{R}_{s;\lambda} (A)$.
Let $\Phi_\zeta$ be the parameter in Chumakin's formula for the generalized resolvents
(see, e.g.~\cite[Theorem 2.7]{cit_46000_Z}), corresponding to
$\mathbf{R}_{u;\zeta} (U_z)$.
These parameters are related by the following formula (see considerations after the proof of Proposition~4.20 and the formula~(4.87)
in~\cite[pp. 280-281]{cit_46000_Z}):
\begin{equation}
\label{f2_52}
F(\lambda) = \Phi_{\frac{ \lambda - z }{ \lambda - \overline{z} }},\quad \lambda\in\Pi_z.
\end{equation}
By the substitution of the expression for $\mathbf{R}_{u;\cdot} (U_z)$ from Chumakin's formula into relation~(\ref{f2_48})
and using~(\ref{f2_52}) we obtain relation~(\ref{f2_40}).

\noindent
$\Box$

Let us return to the investigation of the moment problem. The following relation holds:
$$
(\mathbf{R}_z(A) x_k,x_j)_H =
\frac{1}{z^2 + 1}
\left(
(\mathbf{R}_z(A) y_k,y_j)_H - \gamma_{k+N,j} - z \gamma_{k,j}
\right),
$$
\begin{equation}
\label{f2_57}
z\in \mathbb{C}_+\backslash\{ i \},\ 0\leq k,j \leq N-1.
\end{equation}
This can be proved by the same arguments as a similar property~(28) in~\cite[pp. 393-394]{cit_39000_Z}.

Using~(\ref{f2_30}), (\ref{f2_57}), and relation~(\ref{f2_40}) with $z=i$, we obtain that
all solutions of the moment problem~(\ref{f1_5}) have the following form:
$$ M(\lambda) = (m_{k,j} (\lambda))_{k,j=0}^{N-1}, $$
$$
\int_\mathbb{R} \frac{1}{\lambda - z} dm_{k,j} (\lambda) =
\frac{2i}{ (z^2+1)^2 }
\left\{
\left(
\left[
E_H - \frac{z-i}{z+i}
\left(
U_i(A) \oplus F(z)
\right)
\right]^{-1}
y_k, y_j
\right)_H
\right.
$$
\begin{equation}
\label{f2_60}
\left. + \psi_{j,k} (z)\right\},\quad z\in \mathbb{C}_+\backslash\{ i \},
\end{equation}
where
$F(z)$ is a function from $\mathcal{S}_{a;i}(\mathbb{C}_+; \mathcal{N}_i(A),
\mathcal{N}_{-i}(A))$, and
$$
\psi_{j,k} (z) := -\frac{(z+i)}{2i}
\left(
\gamma_{k+N,j+N} + (z-i)\gamma_{k+N,j} + (z^2 - iz + 1) \gamma_{k,j}
\right),
$$
\begin{equation}
\label{f2_62}
0\leq k,j\leq N-1,\ z\in \mathbb{C}.
\end{equation}
Conversely, an arbitrary function $F(z)\in \mathcal{S}_{a;i}(\mathbb{C}_+; \mathcal{N}_{i}(A),
\mathcal{N}_{-i}(A))$
defines by relation~(\ref{f2_60}) a solution $M(x)$ of the moment problem~(\ref{f1_5}).
Moreover, for different functions from
$\mathcal{S}_{a;i}(\mathbb{C}_+; \mathcal{N}_{i}(A),
\mathcal{N}_{-i}(A))$
there correspond different solutions of the moment problem~(\ref{f1_5}).

\begin{dfn}
\label{d2_1}
Choose an arbitrary $a\in \mathbb{N}$ and $X\in \mathbb{C}_{a\times a}$.
By $\mathrm{S} (\mathbb{C}_+; \mathbb{C}_{a\times a}; X)$ we denote a set of all $\mathbb{C}_{a\times a}$-valued analytic functions $G(z)$ in
$\mathbb{C}_+$, such that
$$ G^*(z) G(z) \leq I_a,\qquad \forall z\in \mathbb{C}_+, $$
and
from the validity of the following relations:
$$ \lim_{ \lambda\in \Pi_i^\varepsilon,\ \lambda\to\infty } G(\lambda) \vec\xi = X \vec\xi, $$
$$ \underline{\lim}_{\lambda\in\Pi_i^\varepsilon,\ \lambda\to\infty}
\left[
|\lambda| (\| \vec\xi \|_{\mathbb{C}_a} - \| G(\lambda) \vec\xi \|_{\mathbb{C}_a})
\right] < +\infty,
$$
for some $\varepsilon$: $0<\varepsilon <\frac{\pi}{2}$, and an element $\vec\xi\in \mathbb{C}_a$,
it follows $\vec\xi = 0$.
\end{dfn}

Consider the operator $X_i=X_i(A)$ forbidden with respect to $A$ (see the corresponding definition in Notations).
Denote by $\widetilde X_i$ the matrix of $X_i$ with respect to the bases $\mathfrak{A}'$, $\mathfrak{A}_v'$.
In order to construct $\widetilde X_i$ by the given moments explicitly, we shall need the following proposition.

\begin{prop}
\label{p2_2}
Let $A$ be a closed symmetric operator in a Hilbert space $H$.
Let
\begin{equation}
\label{f2_65}
S_{\infty;\lambda} = P^H_{\mathcal{N}_{\lambda}(A)}|_{ H\ominus \overline{D(A)} },\qquad \lambda\in \mathbb{R}_e.
\end{equation}
The operator $S_{\infty;\lambda}$ is invertible, for all $\lambda\in \mathbb{R}_e$, and
\begin{equation}
\label{f2_66}
X_{\lambda}(A) = S_{ \infty; \overline{\lambda} } S_{\infty;\lambda}^{-1},\qquad \lambda\in \mathbb{R}_e,
\end{equation}
where $X_{\lambda}(A)$ is the forbidden operator with respect to $A$.

If $H$ is finite-dimensional, then the following relations hold:
\begin{equation}
\label{f2_67}
D(S_{\infty;\lambda}^{-1}) = D(X_{\lambda}(A)) = \mathcal{N}_{\lambda}(A),\qquad \lambda\in \mathbb{R}_e.
\end{equation}
\end{prop}
\textbf{Proof.}
Suppose to the contrary that there exists an element $h\in H\ominus \overline{D(A)}$, $h\not= 0$, such that
$S_{\infty;\lambda} h = P^H_{\mathcal{N}_{\lambda}(A)} h = 0$; $\lambda\in \mathbb{R}_e$. Then
$$ h = P^H_{\mathcal{N}_{\lambda}(A)} h + P^H_{\mathcal{M}_{\lambda}(A)} h = P^H_{\mathcal{M}_{\lambda}(A)} h\in \mathcal{M}_{\lambda}(A), $$
and therefore $h\in (H\ominus \overline{D(A)})\cap \mathcal{M}_{\lambda}(A) = \{ 0 \}$, as it follows from~Corollary~3.5
in~\cite[p. 205]{cit_46000_Z}. We obtained a contradiction, and therefore $S_{\infty;\lambda}$ is invertible.

By~(\ref{f1_5_p2_1}), (\ref{f1_6_p2_1}) we get
\begin{equation}
\label{f2_70}
D(X_\lambda(A)) = S_{\infty;\lambda} (H\ominus \overline{D(A)}) = R(S_{\infty;\lambda});
\end{equation}
$$
X_\lambda(A) S_{\infty;\lambda} h = S_{\infty;\overline{\lambda}} h,\qquad \lambda\in \mathbb{R}_e,\ h\in H\ominus \overline{D(A)}.
$$
Then
$$ X_\lambda(A) y = S_{\infty;\overline{\lambda}} S_{\infty;\lambda}^{-1} y,\qquad \lambda\in \mathbb{R}_e,\ y\in R(S_{\infty;\lambda}). $$
Taking into account relation~(\ref{f2_70}) we obtain~(\ref{f2_66}).

Consider the case of a finite-dimensional $H$. Let $g$ be an arbitrary element of $\mathcal{N}_\lambda(A)$ such that
$g\perp R(S_{\infty;\lambda})$;
$\lambda\in \mathbb{R}_e$. Then
$$ 0 = (g, P^H_{\mathcal{N}_\lambda(A)} h)_H = (g,h)_H,\qquad \forall h\in H\ominus D(A). $$
Then $g\in D(A)$.
By Corollary~3.2 in~\cite[p. 204]{cit_46000_Z} we get $g=0$.
Therefore $\mathcal{N}_\lambda(A)\ominus R(S_{\infty;\lambda}) = \{ 0 \}$, and
taking into account relation~(\ref{f2_70}) we get~(\ref{f2_67}).

\noindent
$\Box$

Return to the investigation of the moment problem. Since we assumed that the moment problem is indeterminate, by Theorem~\ref{t2_1}
we get $\mathfrak{A}_\infty'\not= \emptyset$.

Denote by $\mathcal{M}_{S_{\infty;i}}$ ($\mathcal{M}_{S_{\infty;-i}}$) the matrix of the operator $S_{\infty;i}$ ($S_{\infty;-i}$)
with respect to the bases $\mathfrak{A}'$,$\mathfrak{A}_\infty'$ (respectively with respect to the bases $\mathfrak{A}_v'$,$\mathfrak{A}_\infty'$):

\begin{equation}
\label{f2_75}
\mathcal{M}_{S_{\infty;i}} =
\left( ( S_{\infty;i} f_k', u_l')_H \right)_{0\leq l\leq \delta - 1,\ 0\leq k\leq \kappa'-1} =
\left( ( f_k', u_l')_H \right)_{0\leq l\leq \delta - 1,\ 0\leq k\leq \kappa'-1},
\end{equation}
\begin{equation}
\label{f2_79}
\mathcal{M}_{S_{\infty;-i}} =
\left( ( S_{\infty;-i} f_k', v_l')_H \right)_{0\leq l\leq \delta - 1,\ 0\leq k\leq \kappa'-1} =
\left( ( f_k', v_l')_H \right)_{0\leq l\leq \delta - 1,\ 0\leq k\leq \kappa'-1}.
\end{equation}
Then
\begin{equation}
\label{f2_85}
\widetilde X_i = \mathcal{M}_{S_{\infty;-i}}
\mathcal{M}_{S_{\infty;i}}^{-1}.
\end{equation}

Consider a transformation $\mathbf{T}$ which for an arbitrary function $F(z)\in \mathcal{S}_{a;i}(\mathbb{C}_+; \mathcal{N}_{i}(A), \mathcal{N}_{-i}(A))$
put into correspondence the following $\mathbb{C}_{\delta\times\delta}$-valued function $\mathbf{F}(z)$:
\begin{equation}
\label{f2_90}
\mathbf{F}(z) = \mathbf{T} F =
\left(
(F(z) u_j', v_k')_H
\right)_{0\leq k,j\leq \delta-1},\qquad z\in \mathbb{C}_+.
\end{equation}

It is readily checked that the transformation $\mathbf{T}$ is bijective, and it maps $\mathcal{S}_{a;i}(\mathbb{C}_+; \mathcal{N}_{i}(A), \mathcal{N}_{-i}(A))$
on the whole $\mathrm{S} (\mathbb{C}_+; \mathbb{C}_{\delta\times \delta}; \widetilde X_i)$.

Denote by $\mathcal{M}_{1;z}(F)$ the matrix of the operator $E_H - \frac{z-i}{z+i} \left(
U_i(A)\oplus F(z) \right)$ with respect to the basis $\mathfrak{A}\cup \mathfrak{A}'$; $z\in \mathbb{C}_+$,
$F\in \mathcal{S}_{a;i}(\mathbb{C}_+; \mathcal{N}_{i}(A), \mathcal{N}_{-i}(A))$.
Then
$$ \mathcal{M}_{1;z}(F) = \left(
\begin{array}{cc} A_{0;z} & B_{0;z}(F)\\
C_{0;z} & D_{0;z}(F)\end{array}
\right), $$
where
$$ A_{0;z} =
\left(
\left(
\left(
E_H - \frac{z-i}{z+i} \left( U_i(A)\oplus F(z) \right)
\right) u_k, u_j
\right)_H
\right)_{j,k=0}^{\tau-1} $$
\begin{equation}
\label{f2_95}
=
I_\tau - \left( \frac{z-i}{z+i} \right)
\left(
(v_k,u_j)_H
\right)_{j,k=0}^{\tau-1};
\end{equation}
$$ B_{0;z}(F) =
\left(
\left(
\left(
E_H - \frac{z-i}{z+i} \left( U_i(A)\oplus F(z) \right)
\right) u_k', u_j
\right)_H
\right)_{0\leq j\leq\tau-1,\ 0\leq k\leq\delta-1} $$
\begin{equation}
\label{f2_97}
=
- \left( \frac{z-i}{z+i} \right)
\left(
(F(z) u_k',u_j)_H
\right)_{0\leq j\leq\tau-1,\ 0\leq k\leq\delta-1};
\end{equation}
$$ C_{0;z} =
\left(
\left(
\left(
E_H - \frac{z-i}{z+i} \left( U_i(A)\oplus F(z) \right)
\right) u_k, u_j'
\right)_H
\right)_{0\leq j\leq\delta-1,\ 0\leq k\leq\tau-1} $$
\begin{equation}
\label{f2_99}
=
- \left( \frac{z-i}{z+i} \right)
\left(
(v_k,u_j')_H
\right)_{0\leq j\leq\delta-1,\ 0\leq k\leq\tau-1};
\end{equation}
$$ D_{0;z}(F) =
\left(
\left(
\left(
E_H - \frac{z-i}{z+i} \left( U_i(A)\oplus F(z) \right)
\right) u_k', u_j'
\right)_H
\right)_{0\leq j,k\leq\delta-1} $$
\begin{equation}
\label{f2_101}
=
I_\delta - \left( \frac{z-i}{z+i} \right)
\left(
(F(z) u_k',u_j')_H
\right)_{0\leq j,k\leq\delta-1}.
\end{equation}
Set $\mathbf{F}(z) = (\mathbf{F}_{j,k}(z))_{0\leq j,k\leq\delta-1} := \mathbf{T} F$.
Observe that
\begin{equation}
\label{f2_103}
F(z) u_k' = \sum_{l=0}^{\delta-1} \mathbf{F}_{l,k}(z) v_l',\qquad 0\leq k\leq\delta-1.
\end{equation}
Using~(\ref{f2_103}) we obtain that
\begin{equation}
\label{f2_107}
B_{0;z}(F) =
- \left(\frac{z-i}{z+i}\right)
W \mathbf{F}(z),
\end{equation}
where
\begin{equation}
\label{f2_109}
W := \left(
(v_l',u_j)_H
\right)_{0\leq j\leq\tau-1,\ 0\leq l\leq \delta-1};
\end{equation}
and
\begin{equation}
\label{f2_111}
D_{0;z}(F) =
I_\delta - \left(\frac{z-i}{z+i}\right)
T \mathbf{F}(z),
\end{equation}
where
\begin{equation}
\label{f2_115}
T := \left(
(v_l',u_j')_H
\right)_{0\leq j,l\leq \delta-1}.
\end{equation}
The matrix $\mathcal{M}_{1;z}(F)$ is invertible, since it is the matrix of an invertible operator.
The matrix $A_{0;z}$ is invertible, since it is the matrix of an invertible operator
$$ \left. P^{H}_{ \mathcal{M}_i(A) } \left( E_H - \frac{z-i}{z+i} \left( U_i(A)\oplus F(z) \right) \right) \right|_{ \mathcal{M}_i(A) }
=
E_{ \mathcal{M}_i(A) } - \frac{z-i}{z+i} P^{H}_{ \mathcal{M}_i(A) } U_i(A). $$

By the Frobenius formula
for the inverse of a block matrix we get:
\begin{equation}
\label{f2_117}
\mathcal{M}_{1;z}^{-1}(F) =
\left(
\begin{array}{cc} A_{0;z}^{-1} - \left(\frac{z-i}{z+i}\right) A_{0;z}^{-1} W \mathbf{F}(z) H^{-1}(z) C_{0;z} A_{0;z}^{-1} & *\\
* & *\end{array}
\right),
\end{equation}
where
\begin{equation}
\label{f2_119}
H(z) = I_\delta - \left(\frac{z-i}{z+i}\right)
\left(
T - C_{0;z} A_{0;z}^{-1} W
\right)
\mathbf{F} (z),\quad z\in \mathbb{C}_+,
\end{equation}
and by $(*)$ we denote those blocks which are not important for us.

Let $\{ u_j \}_{j=0}^{\rho-1}$ be those elements, which were obtained during the orthogonalization of the sequence~(\ref{f2_32}), on the first
$N$ steps (i.e. during the orthogonalization of $\{ y_k \}_{k=0}^{N-1}$). Observe that $\rho\geq 1$. In the opposite case, we would have:
$y_k = 0$, $0\leq k\leq N-1$; and by~(\ref{f2_60}) the moment problem would be determinate. That contradicts to our assumptions.
Set
$$ H_\rho^- = \mathop{\rm Lin}\nolimits \{ y_k \}_{k=0}^{N-1} = \mathop{\rm Lin}\nolimits \{ u_j \}_{j=0}^{\rho-1}. $$

Consider the following operator:
$$ J(F;z) =
\left.
P^H_{H_\rho^-}
\left(
E_H - \left(\frac{z-i}{z+i}\right)
\left(
U_i(A)\oplus F(z)
\right)
\right)^{-1}
\right|_{H_\rho^-}:\ H_\rho^-\mapsto H_\rho^-,\qquad z\in \mathbb{C}_+, $$
where $F(z)\in \mathcal{S}_{a;i}(\mathbb{C}_+; \mathcal{N}_{i}(A), \mathcal{N}_{-i}(A))$.
Set
\begin{equation}
\label{f2_125}
k_z = \det \left( (z+i) A_{0;z} \right) = (z+i)^\tau \det A_{0;z},\qquad z\in \mathbb{C}_+,
\end{equation}
\begin{equation}
\label{f2_130}
\widetilde A_{0;z} = \frac{k_z}{z+i} A_{0;z}^{-1},\qquad z\in \mathbb{C}_+.
\end{equation}
Observe that $k_z$ is a scalar polynomial of $z$, while $\widetilde A_{0;z}$ is a matrix polynomial of $z$.
The matrix standing on the intersection of the first $\rho$ rows and the first $\rho$ columns of $\widetilde A_{0;z}$
we denote by $A_{1;z}$.
The matrix standing in the first $\rho$ rows (in the first $\rho$ columns) of $\widetilde A_{0;z}$
we denote by $A_{2;z}$ (respectively by $A_{3;z}$).

\noindent
It is readily checked that the matrix $\mathcal{M}(J(F;z))$ of $J(F;z)$, $z\in \mathbb{C}_+$, with respect to the basis $\{ u_j \}_{j=0}^{\rho-1}$
is equal to
$$
\mathcal{M}(J(F;z)) := \left(
(J(F;z) u_k, u_j)_H
\right)_{j,k=0}^{\rho-1} =
$$
\begin{equation}
\label{f2_135}
\frac{z+i}{k_z} A_{1;z} -
\frac{(z^2+1)}{k_z^2}
A_{2;z} W \mathbf{F}(z) H^{-1}(z) C_{0;z} A_{3;z},\qquad z\in \mathbb{C}_+.
\end{equation}

Consider the following operator, acting from $\mathbb{C}^N$ to $H_\rho^-$:
$$ \mathcal{K} \sum_{n=0}^{N-1} \alpha_n \vec e_n = \sum_{n=0}^{N-1} \alpha_n y_n,\qquad \alpha_n\in \mathbb{C}, $$
where $\vec e_n = (\delta_{n,0},\delta_{n,1},\ldots,\delta_{n,N-1})\in \mathbb{C}^N$.
Let $K$ be the matrix of $\mathcal{K}$ with respect to the orthonormal bases
$\{ \vec e_n \}_{n=0}^{N-1}$ and $\{ u_j \}_{j=0}^{\rho-1}$:
\begin{equation}
\label{f2_140}
K = \left( \left( \mathcal{K} \vec e_k, u_j \right)_H \right)_{0\leq j\leq \rho-1,\ 0\leq k\leq N-1}
= \left( \left(y_k, u_j \right)_H \right)_{0\leq j\leq \rho-1,\ 0\leq k\leq N-1}.
\end{equation}
The right-hand side of~(\ref{f2_60}) may be written as
$$ \frac{2i}{(z^2+1)^2}
\left\{
\left(
\mathcal{K}^* J(F;z) \mathcal{K} \vec e_k, \vec e_j
\right)_{\mathbb{C}^N}
+ \psi_{j,k}(z)
\right\}. $$
Set
$$ \Psi(z) = \left(
\psi_{j,k}(z)
\right)_{0\leq j,k\leq N-1}.
$$
Rewriting~(\ref{f2_60}) in the matrix form we get:
$$ \int_\mathbb{R} \frac{1}{\lambda - z} dM^T (\lambda) =
\frac{2i}{ (z^2+1)^2 }
\left\{ K^* \mathcal{M}(J(F;z)) K + \Psi(z)
\right\}, $$
\begin{equation}
\label{f2_142}
z\in \mathbb{C}_+\backslash\{ i \},\ F(z)\in \mathcal{S}_{a;i}(\mathbb{C}_+; \mathcal{N}_i(A),
\mathcal{N}_{-i}(A)).
\end{equation}
Set
\begin{equation}
\label{f2_145}
\widehat C =
\left(
(v_k,u_j')_H
\right)_{0\leq j\leq\delta-1,\ 0\leq k\leq\tau-1};
\end{equation}
\begin{equation}
\label{f2_150}
\mathbf{A} (z) = (z+i) K^* A_{1;z} K + k_z \Psi(z),
\end{equation}
\begin{equation}
\label{f2_153}
\mathbf{B} (z) = - (z^2+1) K^* A_{2;z} W,
\end{equation}
\begin{equation}
\label{f2_155}
\mathbf{C} (z) = (-z+i) \left(
k_z T + (z-i) \widehat C \widetilde A_{0;z} W
\right),
\end{equation}
\begin{equation}
\label{f2_157}
\mathbf{D} (z) = (z+i) C_{0;z} A_{3;z} K,\qquad z\in \mathbb{C}_+.
\end{equation}
Using these definitions and relation~(\ref{f2_135}) we rewrite~(\ref{f2_142}) in the following form:
$$ \int_\mathbb{R} \frac{1}{\lambda - z} dM^T (\lambda) $$
$$
=
\frac{2i}{ (z^2+1)^2 k_z }
\left\{
\mathbf{A}(z) + \mathbf{B}(z) \mathbf{F}(z)
\left(
(z+i) k_z I_\delta + \mathbf{C}(z) \mathbf{F}(z)
\right)^{-1}
\mathbf{D}(z)
\right\}, $$
\begin{equation}
\label{f2_160}
z\in \mathbb{C}_+\backslash\{ i \},\quad \mathbf{F}=\mathbf{T}F,\ F(z)\in \mathcal{S}_{a;i}(\mathbb{C}_+; \mathcal{N}_i(A),
\mathcal{N}_{-i}(A)).
\end{equation}

\begin{thm}
\label{t2_2}
Let the indeterminate truncated matrix Hamburger moment problem~(\ref{f1_5}) be given.
Let the operator $A$ in the Hilbert space $H$ be constructed as in~(\ref{f2_15}), and
$\mathbf{A}(z)$, $\mathbf{B}(z)$, $\mathbf{C}(z)$, $\mathbf{D}(z)$ be matrix polynomials constructed by~(\ref{f2_150}), (\ref{f2_153}), (\ref{f2_155})
(\ref{f2_157}). Let $k_z$ be the scalar polynomial defined by~(\ref{f2_125}), and $\widetilde X_i$ be the matrix of
the forbidden operator $X_i(A)$ with respect to the bases $\mathfrak{A}'$, $\mathfrak{A}_v'$.
All solutions $M(\lambda)$ of the moment problem~(\ref{f1_5}) have the following form:
$$ M(\lambda) = (m_{k,j} (\lambda))_{k,j=0}^{N-1},\quad  \int_\mathbb{R} \frac{1}{\lambda - z} dM^T (\lambda) $$
$$
=
\frac{2i}{ (z^2+1)^2 k_z }
\left\{
\mathbf{A}(z) + \mathbf{B}(z) \mathbf{F}(z)
\left(
(z+i) k_z I_\delta + \mathbf{C}(z) \mathbf{F}(z)
\right)^{-1}
\mathbf{D}(z)
\right\}, $$
\begin{equation}
\label{f2_160_1}
z\in \mathbb{C}_+\backslash\{ i \},
\end{equation}
where $\mathbf{F}(z)$ is a matrix-valued function
from $\mathrm{S} (\mathbb{C}_+; \mathbb{C}_{\delta\times \delta}; \widetilde X_i)$.
Conversely, an arbitrary matrix-valued function from $\mathrm{S} (\mathbb{C}_+; \mathbb{C}_{\delta\times \delta}; \widetilde X_i)$
defines by~(\ref{f2_160_1}) a solution $M(\lambda)$ of the moment problem~(\ref{f1_5}).
Moreover, for different matrix-valued functions from $\mathrm{S} (\mathbb{C}_+; \mathbb{C}_{\delta\times \delta}; \widetilde X_i)$
by~(\ref{f2_160_1}) there correspond different solutions of the moment problem~(\ref{f1_5}).
\end{thm}
\textbf{Proof.}
The proof follows directly from the previous considerations.

$\Box$

\begin{example}
\label{e2_2}
Consider the moment problem from Example~\ref{e2_1}.
The operator $A$ is defined on $\mathop{\rm Lin}\nolimits\{ x_0, x_1 \}$, and $Ax_0=x_2$, $Ax_1=x_3$.
Then $y_0 = x_2 - ix_0$, $y_1 = x_3 - ix_1$.
Applying the Gram-Schmidt orthogonalization procedure to $y_0,y_1,x_0,x_1$ we obtain that
$u_0 = \frac{\sqrt{3}}{2} (x_2 - ix_0)$, $u_1 = \frac{1}{\sqrt{2}} (x_3 - ix_1)$, and
$u_0' = \left( 3 + \frac{3}{2}i \right) x_0 -
\left(  \frac{3}{2} + i \right) x_2$; $\delta = 1$.
Then $v_0 = \frac{\sqrt{3}}{2} (x_2 + ix_0)$, $v_1 = \frac{1}{\sqrt{2}} (x_3 + ix_1)$.
Applying the Gram-Schmidt orthogonalization procedure to $v_0,v_1,x_0,x_1$ we obtain that
$u_0' = \left( 3 - \frac{3}{2}i \right) x_0 +
\left(  -\frac{3}{2} + i \right) x_2$.
We calculate
$$ W = \frac{\sqrt{3}}{4} i \left(
\begin{array}{cc} 1\\
0\end{array}
\right),\quad T = \frac{1}{2} - \frac{3}{4} i,\quad
K = 2 \left(
\begin{array}{cc} \frac{1}{\sqrt{3}} & 0\\
0 & \frac{1}{\sqrt{2}}\end{array}
\right),
 $$
$$ \Psi (z) = -\frac{(z+i)}{2i} \left(
\begin{array}{cc} \frac{1}{3} z^2 + \left( \frac{1}{2} - \frac{1}{3}i \right) z + \frac{4}{3} - \frac{1}{2} i & 0\\
0 & z^2 +(1-i)z + 2-i\end{array}
\right),\ z\in \mathbb{C}_+; $$
$$ \widehat C = \frac{\sqrt{3}}{4} i (1,0),\quad C_{0;z} =
- \left( \frac{z-i}{z+i} \right) \frac{\sqrt{3}}{4} i (1,0), $$
$$ A_{0;z} = \left(
\begin{array}{cc} 1 - \left( \frac{1}{2} + \frac{3}{4} i \right) \left( \frac{z-i}{z+i} \right) & 0\\
0 &  1 - i \left( \frac{z-i}{z+i} \right)\end{array}
\right),\qquad z\in \mathbb{C}_+; $$
$$ \widetilde A_{0;z} = A_{1;z} = A_{2;z} = A_{3;z} = \left(
\begin{array}{cc} (1-i)(z-1) & 0\\
0 &  \left( \frac{1}{2} - \frac{3}{4} i \right) z - \frac{3}{4} + \frac{3}{2} i \end{array}
\right),\ z\in \mathbb{C}_+; $$
$$ k_z = (1-i) \left(
\left( \frac{1}{2} - \frac{3}{4} i \right)z - \frac{3}{4} + \frac{3}{2} i
\right)
(z-1), \ z\in \mathbb{C}_+. $$
Then
$$ \mathbf{A}(z) = \frac{(z+i)}{2i} \left(
\begin{array}{cc} a_{1,1}(z) & 0\\
0 & a_{2,2}(z)\end{array}
\right),\qquad z\in \mathbb{C}_+, $$
where
$$ a_{1,1}(z) = \frac{1}{3}(1-i)(z-1)\left\{
8i -
\left(
\left( \frac{1}{2} - \frac{3}{4} i \right)z - \frac{3}{4} + \frac{3}{2} i
\right)
\left(
z^2 + \left( \frac{3}{2} - i \right) z + 4 - \frac{3}{2} i
\right)
\right\}, $$
$$ a_{2,2}(z) =
\left(
\left( \frac{1}{2} - \frac{3}{4} i \right)z - \frac{3}{4} + \frac{3}{2} i
\right)
\left(
4i + (-1+i)(z-1)(z^2 + (1-i)z + 2-i)
\right); $$
$$ \mathbf{B}(z) = -\frac{1}{2} (1+i)(z^2+1)(z-1)
\left(
\begin{array}{cc} 1\\
0\end{array}
\right),\quad z\in \mathbb{C}_+; $$
$$ \mathbf{C}(z) =
(1-i)(z-1)(-z+i)
\left(
\left( -\frac{1}{2} - \frac{3}{4} i \right)z + \frac{3}{4} + \frac{3}{2} i
\right),\quad z\in \mathbb{C}_+; $$
$$ \mathbf{D}(z) = -\frac{1}{2} (1+i)(z-1)(z-i) (1,0),\quad z\in \mathbb{C}_+. $$
By the substitution of the latter expressions into relation~(\ref{f2_160}) and after simplifications we get:
$$ \int_\mathbb{R} \frac{1}{\lambda - z} dM^T (\lambda) =
\left(
\begin{array}{cc}
\frac{ -\frac{1}{3} \left(  z - \frac{1}{26} (9-32i)  \right)   }{ (z+i) \left(  z - \frac{3}{13} (8-i)  \right)  } & 0\\
0 & \frac{1}{1-z}\end{array}
\right)
$$
$$ +
\frac{1}{
(z+i)
\left(
\left( \frac{1}{2} - \frac{3}{4} i \right)z - \frac{3}{4} + \frac{3}{2} i
\right)
}
\mathbf{F}(z)
\left(
(z+i)
\left(
\left( \frac{1}{2} - \frac{3}{4} i \right)z - \frac{3}{4} + \frac{3}{2} i
\right)
+
(-z+i)
\right.
$$
$$
\left.
*
\left(
\left( -\frac{1}{2} - \frac{3}{4} i \right)z + \frac{3}{4} + \frac{3}{2} i
\right)
\mathbf{F}(z)
\right)^{-1}
\left(
\begin{array}{cc}
1 & 0\\
0 & 0\end{array}
\right),
 $$
\begin{equation}
\label{f2_165}
z\in \mathbb{C}_+\backslash\{ i \}.
\end{equation}
Let us calculate $\widetilde X_i$. At first, we should calculate $\mathcal{M}_{S_{\infty;\pm i}}$:
$$
\mathcal{M}_{S_{\infty;i}} =
( f_0', u_0')_H = -\frac{3}{4} + \frac{1}{2}i;
$$
$$
\mathcal{M}_{S_{\infty;-i}} =
( f_0', v_0')_H =  -\frac{3}{4} - \frac{1}{2}i.
$$
Then
$$ \widetilde X_i = \mathcal{M}_{S_{\infty;-i}} \mathcal{M}_{S_{\infty;i}}^{-1} =  \frac{5}{13} + \frac{12}{13}i.
$$
Relation~(\ref{f2_165}) establishes a one-to-one correspondence between all functions $\mathbf{F}(z)$ from
$\mathrm{S} \left( \mathbb{C}_+; \mathbb{C}_{1\times 1}; \frac{5}{13} + \frac{12}{13}i \right)$ and all solutions
$M(\lambda) = (m_{k,l}(\lambda))_{k,l=0}^1$ of the moment problem.
In particular, we see that $m_{0,1}(\lambda)=m_{1,0}(\lambda)\equiv 0$,
$m_{1,1}(\lambda) = \left\{
\begin{array}{cc} 0, & \lambda\leq 1\\
1, & \lambda > 1\end{array}
\right.$,
and $m_{0,0}(\lambda)$ can be calculated by the Stieltjes-Perron inversion formula.

\end{example}

\textbf{The determinate case.}
Consider the moment problem~(\ref{f1_5}) with a prescribed sequence $\{ S_n \}_{n=0}^{2d}$ of Hermitian $(N\times N)$ complex matrices,
$d,N\in \mathbb{N}$. Suppose that this moment problem is \textbf{determinate}.
By Theorem~\ref{t2_1} this means that the corresponding operator $A$ in the Hilbert space $H$ is self-adjoint.
If $\mathfrak{A}_\infty = \emptyset$, as it was already mentioned, the unique solution of the moment problem is $M(x)\equiv 0$.
Suppose that $\mathfrak{A}_\infty \not= \emptyset$.
In this case, $\mathfrak{A}_\infty$ is an orthonormal basis in $H$.
By~(\ref{f2_30}) the unique solution can be found from the following relation:
$$ M(\lambda) = (m_{k,j} (\lambda))_{k,j=0}^{N-1},\quad
\int_\mathbb{R} \frac{1}{\lambda - z} dm_{k,j} (\lambda) = ( (A-z E_H)^{-1} x_k, x_j)_H $$
\begin{equation}
\label{f2_170}
= ( \mathcal{R}^* (A- z E_H)^{-1} \mathcal{R} \vec e_k, \vec e_j)_{\mathbb{C}^N},\quad z\in \mathbb{R}_e,
\end{equation}
where the operator $\mathcal{R}$ maps $\mathbb{C}^N$ into $H$:
$$ \mathcal{R} \sum_{n=0}^{N-1} \alpha_n \vec e_n = \sum_{n=0}^{N-1} \alpha_n x_n,\qquad \alpha_n\in \mathbb{C}, $$
and $\vec e_n = (\delta_{n,0},\delta_{n,1},\ldots,\delta_{n,N-1})\in \mathbb{C}^N$.
Let $R$ be the matrix of $\mathcal{R}$ with respect to the orthonormal bases
$\{ \vec e_n \}_{n=0}^{N-1}$ and $\mathfrak{A}_\infty$:
\begin{equation}
\label{f2_175}
R = \left( \left( \mathcal{R} \vec e_k, f_j \right)_H \right)_{0\leq j\leq \kappa-1,\ 0\leq k\leq N-1}
=\left( \left( x_k, f_j \right)_H \right)_{0\leq j\leq \kappa-1,\ 0\leq k\leq N-1}.
\end{equation}
Let $\mathcal{M}_A$ be the matrix of $A$ with respect to the basis
$\mathfrak{A}_\infty$:
\begin{equation}
\label{f2_180}
\mathcal{M}_A = \left( \left( A f_k, f_j \right)_H \right)_{0\leq j,k\leq \kappa-1}.
\end{equation}
Observe that the matrices $R$ and $\mathcal{M}_A$ can be calculated explicitly by the prescribed moments.
Rewrite~(\ref{f2_170}) in the matrix form:
\begin{equation}
\label{f2_185}
\int_\mathbb{R} \frac{1}{\lambda - z} dM^T (\lambda)
= R^* (\mathcal{M}_A - z I_{\kappa})^{-1} R,\quad z\in \mathbb{R}_e.
\end{equation}
The unique solution $M(\lambda)$ of the moment problem can be found from~(\ref{f2_185}) by the Stieltjes-Perron
inversion formula. Moreover, since the solution $M(\lambda)$ has a finite number of points of increase, than
expanding rational functions appearing as entries of the matrix on the right of~(\ref{f2_185}) into simple fractions
we can find $M(\lambda)$.

\section{The truncated matrix Hamburger moment problem with an odd number of moments and having a gap.}

Consider the moment problem~(\ref{f1_5}),(\ref{f1_10}) with a prescribed sequence $\{ S_n \}_{n=0}^{2d}$ of Hermitian $(N\times N)$ complex matrices,
$d,N\in \mathbb{N}$, and a prescribed open set $\Delta$. We shall assume that the corresponding moment problem~(\ref{f1_5}) (i.e. the moment
problem~(\ref{f1_5}) with the same moments) is \textbf{indeterminate}.
In fact, if the corresponding moment problem~(\ref{f1_5}) has no solutions than the  moment problem~(\ref{f1_5}),(\ref{f1_10}) has no
solutions, as well.
If the corresponding moment problem~(\ref{f1_5}) has a unique solution than this solution can be found explicitly, as it was discussed
at the end of the previous section.
Then condition~(\ref{f1_10}) may be verified directly.

Each solution $M(x)$ of the moment problem~(\ref{f1_5}) generates a matrix measure on $\mathfrak{B}(\mathbb{R})$ (i.e.
a $\mathbb{C}_{N\times N}^\geq$-valued function on $\mathfrak{B}(\mathbb{R})$, which is countably additive).
For example, we may define this measure according to~(\ref{f2_25}):
\begin{equation}
\label{f3_5}
M(\delta) = \left(
(\mathbf{E}(\delta) x_k, x_j)_H
\right)_{k,j=0}^{N-1},\qquad \delta\in \mathfrak{B}(\mathbb{R}),
\end{equation}
where $\mathbf{E}(\delta)$ is the spectral measure corresponding to the spectral function $\mathbf{E}_t$ related to the solution $M(x)$.
Any matrix measure $\widetilde M(\delta) = (\widetilde m_{k,j}(\delta))_{k,j=0}^{N-1}$ on $\mathfrak{B}(\mathbb{R})$ satisfying the following relation:
\begin{equation}
\label{f3_7}
\widetilde M((-\infty,x)) = M(x),\qquad x\in \mathbb{R},
\end{equation}
coincides with $M(\delta)$.
In fact, we may consider the following functions:
$$ f_{k,j}(\delta;\alpha;\widetilde M) = (\widetilde M(\delta) (\vec e_k + \alpha\vec e_j), \vec e_k + \alpha\vec e_j )_{\mathbb{C}^N} \geq 0, $$
%$$ = \widetilde m_{k,k}(\delta) + 2\mathop{\rm Re}\nolimits (\alpha \widetilde m_{k,j}(\delta)) + \widetilde m_{j,j}(\delta) $$
%$$ \geq \widetilde m_{k,k}(\delta) - 2\sqrt{ \widetilde m_{k,k}(\delta) } \sqrt{  \widetilde m_{j,j}(\delta) } +
%\widetilde m_{j,j}(\delta) = (\widetilde m_{k,k}(\delta) - \widetilde m_{j,j}(\delta)) \geq 0, $$
where $\alpha\in \mathbb{C}$, $\delta\in \mathfrak{B}(\mathbb{R})$.
The scalar measures $f_{k,j}(\delta;\alpha;\widetilde M)$ and $f_{k,j}(\delta;\alpha; M)$ coincide, since coincide their distribution
functions. On the other hand, the entries of $M(\delta)$ and $\widetilde M(\delta)$ are expressed via $f_{k,j}$ by the
polarization formula. Then $\widetilde M(\delta) = M(\delta)$.
The function on the left of~(\ref{f3_7}) is said to be \textit{the distribution function} of $\widetilde M(\delta)$.

\begin{prop}
\label{p3_1}
Let the indeterminate moment problem~(\ref{f1_5}) with $d\in \mathbb{N}$ be given and the operator $A$ in a Hilbert space $H$ be constructed as
in~(\ref{f2_15}). Let $\Delta\in \mathfrak{B}(\mathbb{R})$ be a fixed set. Let $M(x)$, $x\in \mathbb{R}$, be a solution of the moment
problem~(\ref{f1_5}), $M(\delta),\ \delta\in \mathfrak{B}(\mathbb{R})$, be the matrix measure which is defined by~(\ref{f3_5}) with
the corresponding spectral measure $\mathbf{E}(\delta)$, $\delta\in \mathfrak{B}(\mathbb{R})$.
The following two conditions are equivalent:
\begin{itemize}

\item[{\rm (i)}]
$M(\Delta)=0$;

\item[{\rm (ii)}]
$\mathbf{E}(\Delta)=0$.

\end{itemize}
\end{prop}
\textbf{Proof.}
(ii)$\Rightarrow$(i). It follows directly from relation~(\ref{f3_5}).

\noindent
(i)$\Rightarrow$(ii).
By the construction in~\cite[pp. 282-284]{cit_35000_Z}, the solution $M(x)=(m_{k,j}(\lambda))_{k,j=0}^{N-1}$ is generated by
the left-continuous resolution of unity $\{ \widehat E_\lambda \}_{\lambda\in \mathbb{R}}$ of a self-adjoint operator $\widehat A$
in a Hilbert space $\widehat H\supseteq H$:
$$ m_{k,j}(\lambda) = \left(
P^{\widehat H}_H \widehat E_\lambda x_k, x_j
\right)_H,\qquad 0\leq k,j\leq N-1, $$
and
\begin{equation}
\label{f3_9}
\widehat A = UQU^{-1},
\end{equation}
where $Q$ is an operator of the multiplication by an independent variable in $L^2(M)$, and
$U$ is a unitary operator, mapping $L^2(M)$ on $\widehat H$.

Denote by $\widehat E(\delta)$ ($\delta\in \mathfrak{B}(\mathbb{R})$) the orthogonal spectral measure, corresponding to $\widehat E_t$. Observe that according to
the one-to-one correspondence~(\ref{f2_25}) we have:
$\mathbf{E}_\lambda = P^{\widehat H}_H \widehat E_\lambda$, $\lambda\in \mathbb{R}$, and therefore
\begin{equation}
\label{f3_15}
\mathbf{E}(\delta) = P^{\widehat H}_H \widehat E(\delta),\qquad \delta\in \mathfrak{B}(\mathbb{R}).
\end{equation}%%
By~(\ref{f3_9}) we obtain that
\begin{equation}
\label{f3_20}
\widehat E(\delta) = U E_0(\delta) U^{-1},\qquad \delta\in \mathfrak{B}(\mathbb{R}),
\end{equation}
where $E_0(\delta)$ is the orthogonal spectral measure of $Q$.
Observe that
\begin{equation}
\label{f3_22}
E_0(\delta) f(t) = \chi_\delta(t) f(t),\qquad \delta\in \mathfrak{B}(\mathbb{R}),\ f\in L^2(M),
\end{equation}
where $\chi_\delta(t)$ is the characteristic function of the set $\delta$ (i.e. it is equal to $1$ on $\delta$, and equal to $0$ on $\mathbb{R}
\backslash\delta$).
Since $M(\Delta)=0$, then $\tau_M (\Delta) = \sum_{k=0}^{N-1} m_{k,k} (\Delta) = 0$.
Choose an arbitrary function $f\in L^2(M)$. We may write:
$$ \| \chi_\Delta(t) f(t) \|_{L^2(M)}^2 = \int_{\mathbb{R}} \chi_\Delta(t) f(t) M'_\tau (t) d\tau_M (\chi_\Delta(t) f(t))^* $$
$$ =
\int_{\Delta} f(t) M'_\tau (t) f^*(t) d\tau_M = 0. $$
Therefore $E_0(\Delta)=0$, and by~(\ref{f3_20}),(\ref{f3_15}) we get $\mathbf{E}(\Delta) = 0$.

\noindent
$\Box$

The following two results are simple generalizations of Proposition~4.17 and Theorem~4.21 in~\cite{cit_46000_Z}.

\begin{prop}
\label{p3_2}
Let $A$ be a closed symmetric operator in a Hilbert space $H$, and $\mathbf{E}(\delta)$,
$\delta\in \mathfrak{B}(\mathbb{R})$, be its spectral measure.
Let $\Delta\subseteq \mathbb{R}$ be a fixed open set.
The following two conditions are equivalent:

\begin{itemize}
\item[(i)]  $\mathbf{E}(\Delta) = 0$;
\item[(ii)] The generalized resolvent $\mathbf{R}_z(A)$, corresponding to the spectral measure
$\mathbf{E}(\delta)$, admits an analytic continuation on a set
$\mathbb{R}_e\cup\Delta$.
\end{itemize}
\end{prop}

\begin{thm}
\label{t3_1}
Let $A$ be a closed symmetric operator in a Hilbert space $H$, and
$z\in \mathbb{R}_e$ be an arbitrary fixed point.
Let $\Delta\subseteq \mathbb{R}$ be a fixed open set,
and the following conditions hold:
\begin{equation}
\label{f3_25}
\mbox{the set $\Delta$ consists of points of the regular type of $A$};
\end{equation}
\begin{equation}
\label{f3_27}
P^H_{ \mathcal{M}_{ \overline{z} } (A) } \mathcal{M}_\lambda(A) = \mathcal{M}_{ \overline{z} } (A),\qquad
\forall\lambda\in\Delta.
\end{equation}
Consider an arbitrary  generalized resolvent $\mathbf{R}_{s;\lambda}(A)$ of $A$.
Let $F(\lambda)\in \mathcal{S}_{a;z}(\Pi_{z}; \mathcal{N}_{z}(A),
\mathcal{N}_{\overline{z}}(A))$
corresponds to  $\mathbf{R}_{s;\lambda}(A)$ by the Shtraus formula.
The generalized resolvent $\mathbf{R}_{s;\lambda}(A)$ admits an analytic continuation on
 $\mathbb{R}_e\cup\Delta$ if and only if the following conditions hold:

\begin{itemize}
\item[1)] $F(\lambda)$ admits a continuation on $\Pi_z\cup\Delta$ and this continuation
is continuous in the uniform operator topology;

\item[2)] The continued function $F(\lambda)$ maps isometrically $\mathcal{N}_{z}(A)$ on the whole
$\mathcal{N}_{\overline{z}}(A)$, for all $\lambda\in\Delta$;

\item[3)] The operator $F(\lambda)-\mathcal{W}_\lambda$ is invertible for all $\lambda\in\Delta$, and
$$
(F(\lambda)-\mathcal{W}_\lambda) \mathcal{N}_{z}(A) = \mathcal{N}_{\overline{z}}(A),\qquad
\forall\lambda\in\Delta,
$$
where
\begin{equation}
\label{f3_30}
\mathcal{W}_\lambda P^H_{\mathcal{N}_z(A)} f =
\frac{\lambda - \overline{z}}{\lambda - z} P^H_{\mathcal{N}_{\overline{z}}(A)} f,\qquad
f\in \mathcal{N}_\lambda(A),\ \lambda\in\Delta.
\end{equation}
\end{itemize}
\end{thm}
The proofs of these results follows easily from the fact, that each open subset of $\mathbb{R}$ is a union of at most
countable set of open intervals, and from the above mentioned results in~\cite{cit_46000_Z}.
As it follows from an analogous remark at the end of Section~4 in~\cite{cit_46000_Z}, conditions~(\ref{f3_25}),(\ref{f3_27})
are necessary for the existence of at least one generalized resolvent of $A$,
which admits an analytic continuation on $\mathbb{R}_e\cup\Delta$.

In view of further applications to the moment problem~(\ref{f1_5}), we shall obtain another representation for
the function $\mathcal{W}_\lambda$ from the last theorem. Let $A$,$H$,$\Delta$  satisfy the assumptions of Theorem~\ref{t3_1} with $z=i$.
Consider the following operators:
\begin{equation}
\label{f3_31}
S_\lambda = P^H_{  \mathcal{N}_i(A)  }|_{ \mathcal{N}_\lambda(A)  },\quad
Q_\lambda = P^H_{  \mathcal{N}_{-i}(A)  }|_{ \mathcal{N}_\lambda(A)  },\qquad  \lambda\in\Delta.
\end{equation}
Set $V=U_i(A)$. Observe that
$$ \mathcal{N}_i(A) = N_0(V),\ \mathcal{N}_{-i}(A) = N_\infty(V),\ \mathcal{N}_\lambda(A) = N_{  \frac{\lambda-i}{\lambda+i} }(V),\quad
\lambda\in\Delta. $$
Choose an arbitrary point $\lambda\in\Delta$. Since $\lambda$ is a point of the regular type of $A$, by Proposition~4.18 in~\cite{cit_46000_Z}
with $\lambda_0=i$ we obtain that $\frac{\lambda+i}{\lambda-i}(\in \mathbb{T})$ is a point of the regular type of $V$.
Set $\zeta = \frac{\lambda-i}{\lambda+i}\in \mathbb{T}$. Relation~(\ref{f3_30}), with $z=i$, may be written in the following form:
$$ \mathcal{W}_\lambda P^H_{ N_0(V) } f =
\frac{1}{\zeta} P^H_{ N_\infty(V)} f,\qquad
f\in N_\zeta(V). $$
Comparing this relation with the definition of an operator $W_\zeta$ in~(4.60) in~\cite[p. 270]{cit_46000_Z} we conclude that
$\mathcal{W}_\lambda = W_\zeta$. Moreover, the operators $S$ and $Q$, defined afterwards in~\cite[p. 270]{cit_46000_Z}, coincide
with operators $S_\lambda$ and $Q_\lambda$ from~(\ref{f3_31}), respectively.
By~(4.61) in~\cite[p. 271]{cit_46000_Z} we obtain that $D(\mathcal{W}_\lambda) = \mathcal{N}_i(A)$, $R(\mathcal{W}_\lambda) = \mathcal{N}_{-i}(A)$ and
\begin{equation}
\label{f3_34}
\mathcal{W}_\lambda = \frac{\lambda+i}{\lambda-i} Q_\lambda S_\lambda^{-1},\qquad
\lambda\in\Delta.
\end{equation}
Moreover, the operator $S_\lambda^{-1}$ is bounded and defined on the whole $\mathcal{N}_i(A)$.

In the case of a finite-dimensional Hilbert space, the last theorem can be made more simple. We shall need the following lemma.
\begin{lem}
\label{l3_1}
Let $A$ be a closed symmetric operator in a finite-dimensional Hilbert space $H$, and
$z\in \mathbb{R}_e$ be an arbitrary fixed point.
Let $\Delta\subseteq \mathbb{R}$ be a fixed open set and condition~(\ref{f3_25}) holds.
Then condition~(\ref{f3_27}) holds true.
\end{lem}
\textbf{Proof.}
Observe that
$$ M_\infty(U_z(A)) = \mathcal{M}_{\overline{z}} (A), $$
$$ M_{ \frac{\lambda - z}{\lambda- \overline{z}}  }(U_z(A)) = \mathcal{M}_{\lambda} (A),\qquad \lambda\in\Delta. $$

Denote $V=U_z(A)$. Condition~(\ref{f3_27}) can be written in the following form:
\begin{equation}
\label{f3_35}
P^H_{M_\infty(V)} M_{  \frac{\lambda - z}{\lambda - \overline{z}} } (V) = M_\infty(V),\qquad  \lambda\in\Delta.
\end{equation}
Choose an arbitrary point $\lambda\in\Delta$.
Set $\zeta = \frac{\lambda - z}{\lambda - \overline{z}}$.
Since $\lambda$ is a point of the regular type of $A$, then by Proposition~4.18 in~\cite[p. 277]{cit_46000_Z}
we conclude that $\frac{1}{\zeta}$ is a point of the regular type of $V$.
Since $H$ is finite-dimensional, by Corollary~4.7 in~\cite[p. 268]{cit_46000_Z} we get
$$ (H\ominus R(V)) \dotplus M_\zeta = H. $$
Applying $P^H_{M_\infty(V)}$ to the both sides of the latter equality we obtain relation~(\ref{f3_35}).

\noindent
$\Box$

\begin{thm}
\label{t3_2}
Let $A$ be a closed symmetric operator in a finite-dimensional Hilbert space $H$, and
$z\in \mathbb{R}_e$ be an arbitrary fixed point.
Let $\Delta\subseteq \mathbb{R}$ be a fixed open set,
and condition~(\ref{f3_25}) holds.
Consider an arbitrary  generalized resolvent $\mathbf{R}_{s;\lambda}(A)$ of $A$.
Let $F(\lambda)\in \mathcal{S}_{a;z}(\Pi_{z}; \mathcal{N}_{z}(A),
\mathcal{N}_{\overline{z}}(A))$
corresponds to  $\mathbf{R}_{s;\lambda}(A)$ by the Shtraus formula.
The generalized resolvent $\mathbf{R}_{s;\lambda}(A)$ admits an analytic continuation on
 $\mathbb{R}_e\cup\Delta$ if and only if the following conditions hold:

\begin{itemize}
\item[1)] $F(\lambda)$ admits a continuation on $\Pi_z\cup\Delta$ and this continuation
is continuous in the uniform operator topology;

\item[2)] The continued function $F(\lambda)$ has isometric operators from $\mathcal{N}_{z}(A)$ to
$\mathcal{N}_{\overline{z}}(A)$ as values, for all $\lambda\in\Delta$;

\item[3)] The operator $F(\lambda)-\mathcal{W}_\lambda$ is invertible for all $\lambda\in\Delta$,
where $\mathcal{W}_\lambda$ is from~(\ref{f3_30}).
\end{itemize}
\end{thm}
\textbf{Proof.}
By Lemma~\ref{l3_1} we see that the assumptions of Theorem~\ref{t3_1} are satisfied with $A$,$H$,$\Delta$ and $z$.

\noindent
\textit{Necessity.} If $\mathbf{R}_{s;\lambda}(A)$ admits an analytic continuation on
 $\mathbb{R}_e\cup\Delta$, then conditions~1)-3) of Theorem~\ref{t3_1} are satisfied and they imply conditions~1)-3) of the present
theorem.

\noindent
\textit{Sufficiency.} Let conditions~1)-3) of the present theorem be satisfied.
Choose an arbitrary $\lambda\in\Delta$.
Since $F(\lambda)$ is invertible, the dimension of its range is equal to the dimension of its domain, i.e. to the defect number of $A$.
Then $R(F(\lambda))= \mathcal{N}_{ \overline{z} }(A)$, and condition~2) of Theorem~\ref{t3_1} holds.
A similar argument implies that
$R(F(\lambda)-\mathcal{W}_\lambda)= \mathcal{N}_{\overline{z}}(A)$, and condition~3) of Theorem~\ref{t3_1} holds, as well.
It remains to apply~Theorem~\ref{t3_1}.

\noindent
$\Box$

\begin{prop}
\label{p3_3}
Let the indeterminate moment problem~(\ref{f1_5}) with $d\in \mathbb{N}$ be given and the operator $A$ in a Hilbert space $H$ be constructed as
in~(\ref{f2_15}). Let $\Delta\subseteq\mathbb{R}$ be a fixed open set.
Let $\mathfrak{A}^\lambda = \{ g_j \}_{j=0}^{\widetilde\tau-1}$, $0\leq \widetilde\tau\leq dN$, be an orthonormal basis in $\mathcal{M}_\lambda(A)$,
obtained by the Gram-Schmidt orthogonalization procedure from the following sequence:
$$ x_N - \lambda x_0, x_{N+1} - \lambda x_1,..., x_{dN+N-1} - \lambda x_{dN-1}. $$
Here $\lambda\in \Delta$.
The case $\widetilde\tau=0$ means that $\mathfrak{A}^\lambda = \emptyset$, and $\mathcal{M}_\lambda(A) = \{ 0 \}$.
Then the following conditions are equivalent:

\begin{itemize}

\item[{\rm (a)}]
$\Delta$ consists of points of the regular type of $A$;

\item[{\rm (b)}]
$\mathcal{M}_\lambda(A) \not= \{ 0 \}$, and the matrix $\mathcal{M}_{A-\lambda E_H}$ is invertible, for all $\lambda\in\Delta$.
Here we denote by $\mathcal{M}_{A-\lambda E_H}$ the matrix of the operator $A-\lambda E_H$
with respect to the bases $\mathfrak{A}_\infty$, $\mathfrak{A}^\lambda$.
\end{itemize}

Conditions (a), (b) are necessary for the existence of a solution $M(x)$, $x\in \mathbb{R}$,  of the moment
problem~(\ref{f1_5}),  such that $M(\Delta)=0$, where $M(\delta),\ \delta\in \mathfrak{B}(\mathbb{R})$, is the matrix measure corresponding to the solution.
If conditions~(a),(b) are satisfied then $\widetilde\tau = \tau\geq 1$.
\end{prop}
\textbf{Proof.}
The implication $(b)\Rightarrow (a)$ is obvious.
Conversely, suppose that $\Delta$ consists of points of the regular type of $A$. If $\mathcal{M}_\lambda(A) = \{ 0 \}$, we would get
$D(A)=\{ 0 \}$, $\mathfrak{A}_\infty=\emptyset$, and in this case the moment problem would be determinate (as it was mentioned by the construction
of $\mathfrak{A}_\infty$). Therefore $\mathcal{M}_\lambda(A) \not= \{ 0 \}$. The rest is obvious.

Suppose that there exists a solution $M(x)$, $x\in \mathbb{R}$,  of the moment
problem~(\ref{f1_5}),  such that $M(\Delta)=0$, where $M(\delta),\ \delta\in \mathfrak{B}(\mathbb{R})$, is the corresponding matrix measure.
By Proposition~\ref{p3_1} we get $\mathbf{E}(\Delta)=0$, where $\mathbf{E}(\delta)$, $\delta\in \mathfrak{B}(\mathbb{R})$, is
the corresponding spectral measure.
By Proposition~\ref{p3_2} this means that the corresponding generalized resolvent $\mathbf{R}_z(A)$ admits an analytic continuation on a set
$\mathbb{R}_e\cup\Delta$.
In this case, as it was noticed after Theorem~\ref{t3_1} relation~(\ref{f3_25}) holds.

\noindent
$\Box$

Consider the indeterminate moment problem~(\ref{f1_5}), such as in Proposition~\ref{p3_3}, and suppose that condition~(b) of Proposition~\ref{p3_3}
is satisfied.
We denote by $\mathcal{S}_{a;i}(\Pi_{i}; \mathcal{N}_{i}(A),\mathcal{N}_{-i}(A);\Delta; \mathcal{W})$ a set of
all functions $F(\lambda)$ from $\mathcal{S}_{a;i}(\Pi_{i}; \mathcal{N}_{i}(A), \mathcal{N}_{-i}(A))$,
which satisfy conditions~1)-3) of Theorem~\ref{t3_2} with $z=i$.

Notice that $M_1 := \mathop{\rm Lin}\nolimits\{  g_0, g_1,..., g_{\tau-1}, x_0, x_1,...,x_{N-1} \} = H$.
In fact, by the induction argument we may check that
$\{ x_n \}_{n=0}^{kN+N-1}\subseteq M_1$, $k=0,1,...,d$.

\noindent
Apply the Gram-Schmidt orthogonalization procedure to the following sequence:
$$ g_0, g_1,..., g_{\tau-1}, x_0, x_1,...,x_{N-1}. $$
Observe that the first $\tau$ elements are already orthonormal.
During the orthogonalization of the rest $N$ elements we shall obtain an orthonormal set $\mathfrak{A}^\lambda_+ = \{ g_j' \}_{j=0}^{\delta-1}$.
Observe that $\mathfrak{A}^\lambda_+$ is an orthonormal basis in $\mathcal{N}_\lambda(A)$, $\lambda\in\Delta$.

Let $\mathcal{M}_{S_\lambda}$ ($\mathcal{M}_{Q_\lambda}$) be the matrix of the operator $S_\lambda$ ($Q_\lambda$)
from~(\ref{f3_31}) with respect to the bases $\mathfrak{A}^\lambda_+$, $\mathfrak{A}'$ (respectively, to the bases $\mathfrak{A}^\lambda_+$, $\mathfrak{A}_v'$):
\begin{equation}
\label{f3_40}
\mathcal{M}_{S_\lambda} =
\left(
(S_\lambda g_k', u_j')_H
\right)_{j,k=0}^{\delta-1}
=
\left(
(g_k', u_j')_H
\right)_{j,k=0}^{\delta-1},\qquad  \lambda\in\Delta;
\end{equation}
\begin{equation}
\label{f3_45}
\mathcal{M}_{Q_\lambda} =
\left(
(Q_\lambda g_k', v_j')_H
\right)_{j,k=0}^{\delta-1}
=
\left(
(g_k', v_j')_H
\right)_{j,k=0}^{\delta-1},\qquad  \lambda\in\Delta.
\end{equation}
Denote by $\widetilde{ \mathcal{W} }_\lambda$ the matrix of the operator $\mathcal{W}_\lambda$ from~(\ref{f3_30}) with respect to the
bases $\mathfrak{A}'$, $\mathfrak{A}_v'$. By~(\ref{f3_34}) we get:
\begin{equation}
\label{f3_50}
\widetilde{ \mathcal{W} }_\lambda = \frac{\lambda+i}{\lambda-i}
\mathcal{M}_{Q_\lambda} \mathcal{M}_{S_\lambda}^{-1},\qquad  \lambda\in\Delta.
\end{equation}

\begin{dfn}
\label{d3_1}
Choose an arbitrary $a\in \mathbb{N}$, $X\in \mathbb{C}_{a\times a}$, $\Delta\subseteq \mathbb{R}$, and let
$Y(\lambda)$ be an arbitrary $\mathbb{C}_{a\times a}$-valued function, $\lambda\in\Delta$.
By $\mathrm{S} (\mathbb{C}_+; \mathbb{C}_{a\times a}; X; \Delta; Y)$ we denote a set of all $\mathbb{C}_{a\times a}$-valued functions $G(z)$
from $\mathrm{S} (\mathbb{C}_+; \mathbb{C}_{a\times a}; X)$ which satisfy the following conditions:
\begin{itemize}

\item[{\rm A)}]
$G(z)$ admits a continuation on $\mathbb{C}_+\cup\Delta$, and the continued function $G(z)$ is continuous (i.e. each entry of
$G(z)$ is continuous);

\item[{\rm B)}]
$G^*(z) G(z) = I_a$, for all $z\in\Delta$;

\item[{\rm C)}]
The matrix $G(z) - Y(z)$ is invertible for all $z\in\Delta$.
\end{itemize}
\end{dfn}

It is straightforward to check that the transformation $\mathbf{T}$, defined by~(\ref{f2_90}), maps
$\mathcal{S}_{a;i}(\mathbb{C}_+; \mathcal{N}_{i}(A),\mathcal{N}_{-i}(A);\Delta; \mathcal{W})$
on
$\mathrm{S} (\mathbb{C}_+; \mathbb{C}_{\delta\times \delta}; \widetilde X_i; \Delta; \widetilde{ \mathcal{W} })$.
Here we mean that $\mathbf{T} \emptyset = \emptyset$.

\begin{thm}
\label{t3_3}
Let the indeterminate moment problem~(\ref{f1_5}) with $d\in \mathbb{N}$ be given and the operator $A$ in a Hilbert space $H$ be constructed as
in~(\ref{f2_15}). Let $\Delta\subseteq\mathbb{R}$ be a fixed open set.
Let
condition~(b) of Proposition~\ref{p3_3} be satisfied.
There exists a solution $M(x)$, $x\in \mathbb{R}$,  of the moment
problem~(\ref{f1_5}),  such that $M(\Delta)=0$, where $M(\delta),\ \delta\in \mathfrak{B}(\mathbb{R})$, is the matrix measure corresponding to the solution,
if and only if $\mathrm{S} (\mathbb{C}_+; \mathbb{C}_{\delta\times \delta}; \widetilde X_i; \Delta; \widetilde{ \mathcal{W} })\not=
\emptyset$.
\end{thm}
\textbf{Proof.}
Observe that by Proposition~\ref{p3_3} it follows that $\Delta$ consists of points of the regular type of $A$.

\noindent
\textit{Necessity.}
As in the proof of Proposition~\ref{p3_3} we obtain that the corresponding to the solution $M$ 
generalized resolvent $\mathbf{R}_z(A)$ admits an analytic continuation on a set
$\mathbb{R}_e\cup\Delta$.
Let $F(\lambda)\in \mathcal{S}_{a;i}(\Pi_{i}; \mathcal{N}_{i}(A),\mathcal{N}_{-i}(A))$
corresponds to  $\mathbf{R}_{\lambda}(A)$ by the Shtraus formula.
By Theorem~\ref{t3_2} we conclude that $F(\lambda)\in \mathcal{S}_{a;i}(\Pi_{i}; \mathcal{N}_{i}(A),\mathcal{N}_{-i}(A);\Delta; \mathcal{W})$.
Therefore $\mathrm{S} (\mathbb{C}_+; \mathbb{C}_{\delta\times \delta}; \widetilde X_i; \Delta; \widetilde{ \mathcal{W} })\not=
\emptyset$.

\noindent
\textit{Sufficiency.}
We can choose a function $F(\lambda)\in \mathcal{S}_{a;i}(\Pi_{i}; \mathcal{N}_{i}(A),\mathcal{N}_{-i}(A);\Delta; \mathcal{W})$.
Let $\mathbf{R}_{\lambda}(A) = \mathbf{R}_{s;\lambda}(A)$ be the generalized resolvent of $A$ corresponding to $F$ by the Shtraus formula.
By Theorem~\ref{t3_2}, with $z=i$, we conclude that
$\mathbf{R}_z(A)$ admits an analytic continuation on a set $\mathbb{R}_e\cup\Delta$.
By Proposition~\ref{p3_2} we get $\mathbf{E}(\Delta)=0$. Finally, by
Proposition~\ref{p3_1} we conclude that $M(\Delta)=0$.

\noindent
$\Box$

\begin{thm}
\label{t3_4}
In conditions of Theorem~\ref{t3_3} suppose that
$\mathrm{S} (\mathbb{C}_+; \mathbb{C}_{\delta\times \delta}; \widetilde X_i; \Delta; \widetilde{ \mathcal{W} })\not=
\emptyset$.
All solutions $M(\lambda)$ of the moment problem~(\ref{f1_5}),(\ref{f1_10}) have the form~(\ref{f2_160}) where $\mathbf{F}(z)$ is a matrix-valued function
from 
$\mathrm{S} (\mathbb{C}_+; \mathbb{C}_{\delta\times \delta}; \widetilde X_i; \Delta; \widetilde{ \mathcal{W} })$.
Conversely, an arbitrary matrix-valued function from 
$\mathrm{S} (\mathbb{C}_+; \mathbb{C}_{\delta\times \delta}; \widetilde X_i; \Delta; \widetilde{ \mathcal{W} })$
defines by~(\ref{f2_160}) a solution $M(\lambda)$ of the moment problem~(\ref{f1_5}),(\ref{f1_10}).
Moreover, for different matrix-valued functions from 
$\mathrm{S} (\mathbb{C}_+; \mathbb{C}_{\delta\times \delta}; \widetilde X_i; \Delta; \widetilde{ \mathcal{W} })$
by~(\ref{f2_160}) there correspond different solutions of the moment problem~(\ref{f1_5}),(\ref{f1_10}).
\end{thm}
\textbf{Proof.}
Notice that by Proposition~\ref{p3_3} it follows that $\Delta$ consists of points of the regular type of $A$.

Choose an arbitrary solution $M(\lambda)$ of the moment problem~(\ref{f1_5}),(\ref{f1_10}).
By the construction of formula~(\ref{f2_160}), the corresponding to $M$ matrix-valued function $\mathbf{F}(z)$ is equal to $\mathbf{T}F(z)$,
where $F(z)\in \mathcal{S}_{a;i}(\Pi_{i}; \mathcal{N}_{i}(A),\mathcal{N}_{-i}(A))$ corresponds, 
by the Shtraus formula, to the generalized resolvent $\mathbf{R}_z(A)$, related to $M$.
Repeating the arguments in the proof of the necessity of Theorem~\ref{t3_3} we obtain that
$F(\lambda)\in \mathcal{S}_{a;i}(\Pi_{i}; \mathcal{N}_{i}(A),\mathcal{N}_{-i}(A);\Delta; \mathcal{W})$.
Then $\mathbf{F}\in \mathrm{S} (\mathbb{C}_+; \mathbb{C}_{\delta\times \delta}; \widetilde X_i; \Delta; \widetilde{ \mathcal{W} })$.

On the other hand, choose an arbitrary function
$\mathbf{F}\in \mathrm{S} (\mathbb{C}_+; \mathbb{C}_{\delta\times \delta}; \widetilde X_i; \Delta; \widetilde{ \mathcal{W} })$.
Set $F = \mathbf{T}^{-1} \mathbf{F}\in \mathcal{S}_{a;i}(\mathbb{C}_+; \mathcal{N}_{i}(A),\mathcal{N}_{-i}(A);\Delta; \mathcal{W})$.
Repeating the arguments in the proof of the sufficiency of Theorem~\ref{t3_3} we obtain that
$M(\Delta)=0$, where $M$ is the corresponding to $\mathbf{F}$ solution of the moment problem.

The fact that the correspondence between solutions and
matrix-valued functions from $\mathrm{S} (\mathbb{C}_+; \mathbb{C}_{\delta\times \delta}; \widetilde X_i; \Delta; \widetilde{ \mathcal{W} })$
is one-to-one follows from Theorem~\ref{t2_2}.

\noindent
$\Box$

\begin{example}
\label{e3_1}
Consider the moment problem from Example~\ref{e2_1}.
Set $\Delta =(-1,1)$, and consider the corresponding moment problem~(\ref{f1_5}),(\ref{f1_10}) (i.e. the moment problem
from Example~\ref{e2_1} with an additional constraint~(\ref{f1_10})).
Choose an arbitrary $\lambda\in(-1,1)$.
Apply the Gram-Schmidt orthogonalization procedure to the following sequence:
$$ x_2 - \lambda x_0, x_3 - \lambda x_1. $$
We get
$$ g_0 = \frac{ \sqrt{3} }{ \sqrt{ \lambda^2 - 3\lambda + 3 } } (x_2 - \lambda x_0), $$
$$ g_1 = \frac{1}{1-\lambda} (x_3 - \lambda x_1). $$
Thus, $\mathfrak{A}^\lambda = \{ g_0, g_1 \}$, $\widetilde\tau = 2$, is an orthonormal basis in $\mathcal{M}_\lambda(A)$.
Then
$$ \mathcal{M}_{A-\lambda E_H} =
\left(
((A-\lambda E_H) f_k, g_j)_H
\right)_{0\leq j,k\leq 1}
=
\left(
\begin{array}{cc} \sqrt{ \lambda^2 - 3\lambda + 3 } & 0\\
0 & 1-\lambda\end{array}
\right). $$
Thus, condition~(b) of Proposition~\ref{p3_3} is satisfied.
Apply the Gram-Schmidt orthogonalization procedure to the following sequence:
$$ g_0, g_1, x_0, x_1. $$
We get $g_0' = \frac{1}{ \sqrt{ \lambda^2 - 3\lambda + 3 }   }
\left(
(-3\lambda +6) x_0 + (2\lambda - 3) x_2
\right)$.
Then
$$ \mathcal{M}_{S_\lambda} =
\frac{1}{ 4\sqrt{ \lambda^2 - 3\lambda + 3 }   }
\left(
(-3 +2i) \lambda + 6-3i
\right), $$
$$\mathcal{M}_{Q_\lambda} =
\frac{1}{ 4\sqrt{ \lambda^2 - 3\lambda + 3 }   }
\left(
(-3 -2i) \lambda + 6+3i
\right), $$
and
$\widetilde{ \mathcal{W}_\lambda } =
\frac{ (\lambda+i) ((-3-2i)\lambda + 6 + 3i)  }
{ (\lambda-i) ((-3+2i)\lambda + 6 - 3i)  }$.
Observe that the function $\mathbf{F}(\lambda)=1$, $\lambda\in \mathbb{C}_+$,
belongs to the set
$\mathrm{S} \left( \mathbb{C}_+; \mathbb{C}_{1\times 1}; \frac{5}{13} + \frac{12}{13}i; (-1,1); \frac{ (\lambda+i) ((-3-2i)\lambda + 6 + 3i)  }
{ (\lambda-i) ((-3+2i)\lambda + 6 - 3i)  } \right)$.
Thus, the moment problem~(\ref{f1_5}),(\ref{f1_10}) has a solution.

\noindent
Relation~(\ref{f2_165}) establishes a one-to-one correspondence between all functions $\mathbf{F}(z)$ from
$\mathrm{S} \left( \mathbb{C}_+; \mathbb{C}_{1\times 1}; \frac{5}{13} + \frac{12}{13}i; (-1,1); \frac{ (\lambda+i) ((-3-2i)\lambda + 6 + 3i)  }
{ (\lambda-i) ((-3+2i)\lambda + 6 - 3i)  } \right)$ and all solutions
$M(\lambda)$ of the moment problem.

\end{example}

\vspace{1cm}

Sergey M. Zagorodnyuk

School of Mathematics and Mechanics

Karazin Kharkiv National University

Kharkiv, 61022

Ukraine

e-mail: Sergey.M.Zagorodnyuk@univer.kharkov.ua

\vspace{1cm}
%\newpage

{\normalsize
\begin{center}
\bf
ON THE TRUNCATED MATRIX POWER MOMENT PROBLEM WITH AN OPEN GAP
\end{center}
\begin{center}
\bf
S.M. Zagorodnyuk
\end{center}

In this paper we study the truncated power moment problem with an odd number of prescribed moments. A Nevanlinna-type formula is derived for
this moment problem in the case when the moment problem has more than one solution (the indeterminate case).
The coefficients of the corresponding linear fractional transformation
are expressed explicitly in terms of the given moments.
In the determinate case the solution to the moment problem is constructed, as well.
Those solutions of the moment problem, which satisfy the following condition: $M(\Delta)=0$, where $\Delta$ is a prescribed open subset of
$\mathbb{R}$, and $M$ is the matrix measure generated by a solution of the moment problem, are described.

}

%\vspace{1cm}
MSC 2010: 44A60, 47A57, 30E05

Key words and phrases: moment problem, matrix measure, generalized resolvent.
%18_01_2014_2_05
}

\end{document}